\documentstyle[leqno,amssymb]{bcp}
\font\tenfrak=eufm10
\font\sevenfrak=eufm7
\font\fivefrak=eufm5
        \newfam\frakfam \def\frak{\fam\frakfam\tenfrak} 
                \textfont\frakfam=\tenfrak
\font\tenDDl=msbm10  
\font\sevenDDl=msbm7 
\font\fiveDDl=msbm5 
        \newfam\DDlfam  
                \textfont\DDlfam=\tenDDl 
\scriptfont\DDlfam=\sevenDDl
        \scriptscriptfont\DDlfam=\fiveDDl
\scriptfont\frakfam=\sevenfrak  
        \scriptscriptfont\frakfam=\fivefrak

\long\def\nodo#1{}
\def\genfd{{\bf k}}
\def\mat22#1#2#3#4{ \left(\begin{array}{cc} #1 & #2\\ #3 & #4
\end{array}\right) } 
\newcommand{\id}{{\rm id}}

 \def\Fhd#1#2{\smash{\mathop{\hbox to 14mm{\rightarrowfill}}
\limits^{\scriptstyle#1}_{\scriptstyle#2}}}

\def\Fhg#1#2{\smash{\mathop{\hbox to 14mm{\leftarrowfill}}
\limits^{\scriptstyle#1}_{\scriptstyle#2}}}

 \def\fhd#1#2{\smash{\mathop{\hbox to 8mm{\rightarrowfill}}
\limits^{\scriptstyle#1}_{\scriptstyle#2}}}

\def\fhg#1#2{\smash{\mathop{\hbox to 8mm{\leftarrowfill}}
\limits^{\scriptstyle#1}_{\scriptstyle#2}}}

\def\diagram#1{\def\normalbaselines{\baselineskip=0pt
\lineskip=10pt\lineskiplimit=1pt}   \matrix{#1}}

\def\filtf{{\cal L}} 

\def\LaTeX{L\kern -.36em\raise .3ex\hbox{\sc a}\kern -.15em T\kern -.1667em%
\lower .7ex\hbox{E}\kern -.125em X}

\font\tt=cmtt10
\begin{document}

\keywords{coinvariants,
Hopf algebra, quantum group, quantum principal bundle, localization, 
Ore set, matrix bialgebra, noncommutative Gauss decomposition}
\mathclass{Primary 14A22; Secondary 16W30,14L30,58B32.}
\abbrevauthors{Zoran {\v{S}}koda} 
\abbrevtitle{Localizations for coset spaces}

\title{Localizations for construction \\
of quantum coset spaces}

\author{Zoran { \v{S}}koda}
\address{Department of Mathematics, Indiana University\\
Rawles Hall - Bloomington, IN 47405, USA\\
E-mail: zskoda@indiana.edu}

\maketitlebcp

\abstract{Viewing comodule algebras as the noncommutative
analogues of affine varieties with affine group actions,
we propose rudiments of a localization approach to
nonaffine Hopf algebraic quotients of noncommutative affine varieties
corresponding to comodule algebras. After reviewing basic
background on noncommutative localizations, we introduce
localizations compatible with coactions. Coinvariants
of these localized coactions give local information about
quotients. We define Zariski locally trivial
quantum group algebraic principal and associated bundles. 
Compatible localizations induce localizations 
on the categories of Hopf modules.
Their interplay with the functor of taking coinvariants
and its left adjoint is stressed out. 

Using localization approach,
we constructed a natural class of examples of quantum coset
spaces, related to the quantum flag varieties of type A of other authors.
Noncommutative Gauss decomposition via quasideterminants
reveals a new structure in noncommutative matrix bialgebras.
In the quantum case, calculations with
quantum minors yield the structure theorems.
}

\section*{Notation.} Ground field is $\genfd$ and we assume it
is of characteristic zero.
If we deal just with one $\genfd$-Hopf algebra, say ${\cal B}$, the
comultiplication is 
$\Delta : {\cal B} \rightarrow {\cal B}\otimes {\cal B}$, 
unit map $\eta : \genfd \rightarrow {\cal B}$, 
counit $\epsilon :  {\cal B}\rightarrow \genfd$,
multiplication $\mu :  {\cal B}\otimes {\cal B}\rightarrow {\cal B}$, 
and antipode (coinverse) is $S : {\cal B}\rightarrow {\cal B}$.
Warning: letter $S$ often stands for a generic Ore set.
We use~\cite{ Montg, Majid, KlymikSchmud, Sch:lec}
Sweedler notation $\Delta(h) = \sum h_{(1)} \otimes h_{(2)}$
with or without explicit summation sign, as well as its extension
for coactions: $\rho(v) = \sum v_{(0)} \otimes v_{(1)}$, where
the zero-component is in the comodule and nonzero component(s)
in the coalgebra. 
An entry symbol and name of a matrix will match, 
except for upper vs. lower case, 
e.g. $G = (g^i_j)$; and $G^I_J$ will be a submatrix
with row multilabel $I = (i_1,\ldots, i_k)$ and column multilabel
$J = (j_1,\ldots, j_k)$. As a rule, row labels are placed
as superscripts and column labels as subscripts.

\section*{1. Introduction.} 
There is an antiequivalence~\cite{Jantzen:AG,Sch:lec} between
the category of affine group schemes over $\genfd$ 
and the category of commutative Hopf algebras over $\genfd$.
In the framework of affine group $\genfd$-varieties, the 
corresponding Hopf algebras are constructed as algebras of
regular functions on the underlying group variety.
Hence, taking the viewpoint of 
noncommutative geometry~\cite{Connes:book}
we view Hopf algebras as {\it noncommutative affine group varieties}
~\cite{Drinf:ICM,ManinQGNG, ParshallWang}.   
\vskip .05in
Groups are useful as they describe the notion of symmetries:
they act on spaces.
A {\bf $B$-variety} is an algebraic variety $E$ 
with a regular action $\nu : E \times B \rightarrow E$
of an algebraic group $B$. Hopf algebra ${\cal O}(B)$ of regular
functions on $B$ coacts on the algebra of regular
functions on ${\cal O}(E)$ via coaction 
$\rho: {\cal O}(E) \rightarrow {\cal O}(E) \otimes {\cal O}(B)$
given by formula $\rho(e)(b) = \nu(e,b)$.
This ${\cal O}(B)$-comodule structure on ${\cal O}(E)$ is
compatible with $\genfd$-algebra structure on ${\cal O}(E)$ in the sense
that the coaction is an $\genfd$-algebra map, i.e. ${\cal O}(E)$
is an ${\cal O}(B)$-{\it comodule algebra}. Hence, in the noncommutative
setup, comodule algebras are viewed as noncommutative
$B$-varieties. 
\vskip .05in
When we denote an algebra (Hopf algebra)
by a caligraphic letter, say ${\cal E}$ (or ${\cal B}$) then the
letter still suggests the underlying ``variety'', as in
${\cal E} = {\cal O}(E)$, however we replace ${\cal O}(E)$
by ${\cal E}$ precisely when we allow (though do not prescribe) 
for noncommutative algebras.
\vskip .05in
Any function invariant on orbits of action of a group on a set
can be viewed as a function on the set of orbits. Invariant
functions are coinvariants in the algebra of all functions
with respect to the coaction of Hopf algebra of functions
on the group. This is {\it mutatis mutadis} true in various
setups -- finite sets, topological spaces, 
affine vs. nonaffine algebraic varieties, 
so we did not say it fully precisely. This is very important.
Namely, if we reverse the question and ask could we describe
the space of orbits (quotient space) by coinvariants,
then, already in the commutative situation, the answer
depends much on the category chosen, and even if we
start with a nice-behaved category, the {\it natural quotients
should often be constructed in a larger category}. 
One such phenomenon stems from the fact 
that observables on the quotient space have 
singularities. For example, if $G$ is an affine algebraic
group and $B$ a closed algebraic subgroup, quotient $G/B$
is, in general, not an affine variety, 
but it is always a quasiprojective variety~\cite{Borel}. 
In particular, it is often a projective variety
where the only global regular functions are constants.
However, there are many observables with the
singularity locus of higher codimension and
regular behaviour away from singularities.
Hence it may be sufficient to introduce
regular functions on Zariski open subsets of quotient spaces
(complements of possible singularity loci).
Idea of locally defined quotients is
one of the starting points of
the geometric invariant theory. 
\vskip .015in
This survey gives an overview
of efforts, at this point mainly of the present author,
to access and use the local information on noncommutative
quotients, mainly in the case of noncommutative
coset spaces. The following
are crucial observations in this programme:
\begin{itemize}

\item Noncommutative localizations are used to
replace Zariski open subsets.

\item Already in the commutative case, open sets in the quotient $E/B$ 
correspond (via projection $E \rightarrow E/B$), not
to arbitrary, but only to {\it $B$-invariant} open sets in $E$.
To address this issue, in dual language, we introduce and study 
a notion of {\bf compatibility} 
of a noncommutative localization with coaction. 
Coaction naturally extends to compatible localization.
We obtain the {\bf localized coaction}.
Any compatible localization of
a ${\cal B}$-comodule algebra ${\cal E}$ induces
a localization functor from the category of 
those modules over ${\cal E}$ which are also
${\cal B}$-comodules in a compatible way (relative Hopf modules).
They are analogues of $B$-equivariant (quasicoherent)
sheaves on $E$. 

\item We study {\bf localized coinvariants} i.e. coinvariants
for the localized coaction.

\item In the case of a noncommutative coset spaces, 
the existence of
a covering by compatible localizations with large 
algebras of localized coinvariants justifies calling
the latter charts in a coset space. Large is here in
the sense of ability to perform descent, e.g.
if localization has induced structure
of a faithfully flat Hopf-Galois extension.

\item Noncommutative Gauss decomposition for
matrix bialgebras suggests natural candidates
for (covers by) coaction compatible 
Cohn's, and, in favorable cases, Ore localizations,
such that localized algebras 
can be trivialized as ${\cal B}$-bundles
over the algebras of localized coinvariants. 
This provides a natural class of noncommutative candidates
for coset spaces.
 
\item Calculations with quasideterminants and,
in the case of quantum groups, quantum minors,
are useful techniques to study
the above mentioned examples.
\end{itemize}

\section*{2. Commutative localizations (motivation).}

In (commutative) geometry localization appears as a means to
\begin{itemize}
\item 1. pass from a space to an open subset of the space;
\item 2. pass to a different space reflecting
only infinitesimal neighborhood of a point 
or a subvariety. 
\end{itemize}
In this paper we concentrate on the former flavor of localization.
In the language of algebras, localizing can be done
by introducing functions defined only locally.
For affine varieties that means 
{\it introducing inverses} of those elements in the algebra
whose zero set lies outside of the local set.
There were also some attempts to use a localization
for algebras of continuous functions~\cite{fre-loc}, 
and that type of procedure may be useful 
for extensions of the present work to
operator-algebraic, rather than algebraic setup.
\vskip .01in
The localized ring has a simpler structure than
the original ring. Indeed, if $f$ was a generator of an ideal $I$
different from the whole ring $R$, then having $f^{-1}$ in localized
ring means that $f$ does not generate a proper ideal any more,
as $f^{-1}f = 1$. Hence, {\it localization
kills ideals}. In particular it kills prime and maximal ideals
and, as those correspond to points of schemes and varieties
respectively, it removes some points from the space and the space
gets smaller or "localized". For localization at a point
we obtain a local ring.
\vskip .01in
If we introduce inverses of functions, we still know how to
multiply them: pointwise. Noncommutative algebras
are not algebras of functions on a genuine space consisting
of points only, so we do not have 
apriori fully satisfactory recipe 
for how to multiply the newly introduced inverses
with other elements. There is an important case when such a recipe
is known and elementary. That is the case of inverting
all elements belonging to a given subset $S$ in the ring $R$
of special kind, so called Ore set. This is {\it Ore localization}.

\section*{3. Ore localizations.} ~\cite{ore:31, elizarov, Skoda:locelec}
A semigroup $R$ with unit is called a monoid. 
A subset $S$ of a monoid $R$ is called multiplicative 
if $1 \in R$ and whenever $s_1,s_2 \in S$ then $s_1 s_2 \in S$. 
Let $R$ be a (noncommutative) unital ring. We can also view it as
a monoid with respect to multiplication.
A multiplicative set $S \subset R\backslash \{0\}$ 
is called an {\bf {left Ore set}} 
if the following {\bf left Ore conditions} are satisfied:
\begin{itemize}
\item $(\forall s \in S \mbox{ } \forall r \in R 
\mbox{ }\exists s' \in S \mbox{ }\exists r' \in R) ( r' s = s' r )$
(left Ore condition proper)
\item $(\forall n \in R \mbox{ }\forall s \in S) 
((ns = 0) \Rightarrow \exists s' \in S (s'n = 0))$ (left reversibility)
\end{itemize}
The left reversibility condition can be restated also as
\[ (\forall n_1,n_2 \in R \mbox{ }\forall s \in S) \,\,
((n_1 s = n_2 s) \Rightarrow \exists s' \in S \,(s'n_1 = s'n_2)), \]
what has the advantage that it makes sense for arbitrary monoids,
as well as, once the quantifiers are rewritten
with care to appropriate source and target matching, also for groupoids,
and categories. In the latter case we obtain ``left calculus of fractions''
rather than left Ore set, but the construction of localization (this time
of a category) and accompanying
proofs may proceed essentially the same way as for Ore sets.
\vskip .01in
For a left Ore set $S$ in a monoid $R$ define the monoid
$S^{-1}R$ of {\bf left fractions} as follows. 
As a set, $S^{-1}R := S\times R/\sim$, where 
$\sim$ is the following relation of equivalence:
\[ (s,r) \sim (s',r') \,\,\Leftrightarrow\,\,
(\exists \tilde s\in S \,\,\exists \tilde{r}\in R) \,\,
(\tilde{s}s' = \tilde{r}s\, \mbox{ and } \,\tilde{s}r' = \tilde{r}r).\]
The class of equivalence of $(s,r)$ is denoted $s^{-1}r$ and
called a left fraction. The multiplication is defined by
$s_1^{-1}r_1\cdot s_2^{-1}r_2 = (\tilde{s}s_1)^{-1} (\tilde{r}r_2)$ where 
$\tilde{r} \in R, \tilde{s} \in S$
satisfy $\tilde{r}s_2 = \tilde{s}r_1$ 
(one should think of this, though it is not yet
formally justified at this point, as 
$\tilde{s}^{-1}\tilde{r} = r_1 s_2^{-1}$, what enables to put inverses one
next to another and then the multiplication rule is obvious).
If the monoid $R$ is a ring, then we can extend the addition to
$S^{-1}R$ too. Suppose we are given two fractions
with representatives $(s_1,r_1)$ and $(s_2,r_2)$. Then by the left Ore
condition we find $\tilde{s} \in S$, $\tilde{r}\in R$ such that
$\tilde{s} s_1 = \tilde{r} s_2$. The sum is then defined 
\[ s_1^{-1} r_1 + s_2^{-1} r_2 \,:=\, 
 (\tilde{s}s_1)^{-1} (\tilde{s}r_1 + \tilde{r}r_2) \]
It is a long and at points tricky to work out all the details of this
definition. One has to show that $\sim$ is indeed relation of equivalence,
that the operations are well defined, and that $S^{-1}R$ is indeed
a ring. Even the commutativity of addition needs work.
At the end one shows that $i = i_S : R \rightarrow S^{-1}R$ given by
$i(r) = 1^{-1}r$ is a homomorphism of rings, which is 1-1 iff
the 2-sided ideal $I_S = \{ n \in R \,|\,\exists s \in S,\, sn = 0\}$
is zero.

If $S$ is left Ore, then we call the pair $(i,S^{-1}R)$ the left
Ore localization of $R$ with respect to $S$.
It has a {\bf universal property},
namely, it is a universal object in category ${\cal C} = {\cal C}(R,S)$ 
whose objects are pairs $(j,Y)$, where $j : R \rightarrow Y$ is a map into a
ring $Y$ such that image $j(S)$ of $S$ consists of units, and
the morphisms $\alpha : (j,Y) \rightarrow (j',Y')$
are maps of rings $\alpha : Y \rightarrow Y'$ such that $\alpha \circ j = j'$.
A universal object in ${\cal C}$ may exist when $S$ is not left Ore,
for example when $S$ is right Ore and not left Ore. In fact,
the universal object is a left Ore localization iff it lies in the
full subcategory ${\cal C}^l$ of ${\cal C}$ whose objects $(j,Y)$
satisfy 2 additional conditions: 
$j(S)^{-1}j(R) = \{(j(s))^{-1}j(r)\,|\, s \in S, r\in R\}$ is
a subring in $Y$ and ${\rm ker}\,j = I_S$. Hence
$(i,S^{-1}R)$ is universal in ${\cal C}^l$, and that characterizes
it, but the universality in $\cal C$, although
not characteristic, appears to be more useful in practice.

If $M$ is a left $R$-module then 
$S^{-1}M = S^{-1}R \otimes_R M$ is also a left module. 
This is the recipe for {\bf Ore localization of modules}.
Correspondence $Q_S : M \mapsto S^{-1}M$ is an 
exact endofunctor in category of left $R$-modules, called
the localization functor. For given $M$,
rule $i_{S,M} : M \rightarrow S^{-1}M$ given by $m \mapsto 1 \otimes m$
is an $R$-module map, called the localization map.

\vskip .013in
Ore sets are relatively rare and also hard to single out. 
In practice, Ore condition is checked on a suitable set of generators of
a ring versus a suitable set of generators of the Ore set.
One often uses induction arguments, recursively satisfying
the Ore condition.

\section*{4. Ore vs. Gabriel localizations.} 
\cite{Rosen:88, JaraVerschorenVidal}
This section could be skipped in first reading as only few
remarks in the paper depend on it. The modern viewpoint on localization
as touched upon here is however essential for the current research
in this area.

\vskip .023in
A {\bf lattice} is a poset $(W,\succ)$ such that for any two elements
$z_1,z_2$ the least upper bound $z_1 \vee z_2$
and the greatest lower bound $z_1 \wedge z_2$ exist. 
In other words, binary operations 
{\it meet} $\wedge$ and {\it join} $\vee$ are everywhere defined.
A poset is {\bf bounded} if it contains a maximum and a minimum element,
which we denote $1$ and $0$ respectively.
A {\bf filter} in a bounded lattice $(W,\succ)$
is a subset $\filtf \subset W$ such that
$1\in \filtf$, $0 \notin \filtf$,
$(z_1,z_2 \in \filtf\Rightarrow z_1 \wedge z_2 \in \filtf)$
and $(z \in \filtf, z' \succ z \Rightarrow z' \in \filtf)$.

\vskip .027in
For any subset $w \subset R$, and any left ideal $J$, denote
$(J : w) = \{ z \in R\,|\,zw \in J\}$. It is also a left ideal.
Let $I_l R$ be the preorder category of left ideals in a ring $R$
with respect to inclusion preorder. It is a lattice. 
For the localization questions another partial order $\succ$ on
$I_l R$ is sometimes better. Namely, $K \succ J$ iff either $J \subset K$
or there exist a finite subset $w \subset R$ such that $(J:w) \subset K$.
Any filter in $( I_l R, \succ )$ is called a {\bf uniform filter}.

\vskip .025in
For an Ore set $S \subset R$ consider 
$\filtf_S = \{ J\mbox{ left ideal in } R\,|\, 
J  \cap S \neq 0 \} \subset I_l R.$

\vskip .01in
Left Ore condition implies at once that $\filtf_S$ can
equivalently be defined by
\begin{equation}\label{eq:Gabriel-filter-for-S}
 \filtf_S = \{ J\mbox{ left ideal in } R
\,\,|\, \forall r \,\,(J : r) \cap S \neq 0
   \}. 
\end{equation}
For any multiplicative subset $S \subset R$, not necessarily
left Ore, formula~(\ref{eq:Gabriel-filter-for-S}) defines
a {\bf Gabriel filter} $\filtf_S$ of left ideals in $R$.
It is a uniform filter.

\vskip .027in
To any Gabriel filter $\filtf$, one associates 
endofunctor $\sigma_\filtf$ on the category 
of left $R$-modules by 
\[ \sigma_\filtf(M) = \{ m \in M \,|\, \exists J \in \filtf, \,J m = 0\}.
\]
Equivalently, 
$\sigma_\filtf(M) = {\rm lim}_{J \in \filtf} {\rm Hom}_R (R/J,M)$.
For example, if $\filtf = \filtf_S$ where $S$ is Ore, then
$\sigma_\filtf(R) = I_S$ (see section 3).

\vskip .03in
A {\bf subobject} in a category is an equivalence class of monomorphisms.
Functor $F$ is a {\bf subfunctor} of functor $G$ if
${\rm in}_M : F(M)\hookrightarrow G(M)$ is a subobject
and the inclusions ${\rm in}_M : F(M) \hookrightarrow G(M)$
form a natural transformation of functors ${\rm in}: F \rightarrow G$.
Explicitly, ${\rm in}_N \,F(f)(F(M)) =  G(f)({\rm in}_M F(M))$ 
for $f: M \rightarrow N$.

\vskip .026in
If ${\cal A}$ is any Abelian category, then a
subfunctor $\sigma$ of the identity 
(i.e. $\sigma(M)\subset M$ and $\sigma(f)(\sigma(M)) = f(\sigma(M))$)
with the property 
$\sigma(M/\sigma(M)) = 0$ is called a {\bf preradical} in ${\cal A}$.
A {\bf radical} is a left exact preradical.
It follows that $\sigma_\filtf$ is an {\bf idempotent radical} in
the category of left $R$-modules i.e. it is a radical 
and $\sigma_\filtf(\sigma_\filtf(M)) = \sigma_\filtf(M)$.

\vskip .007in
To any Gabriel filter $\filtf$, one associates a localization
endofunctor $Q_{\filtf}$ on the category of left modules by
the formula
\[ Q_\filtf(M) = 
{\rm lim}_{J \in \filtf} {\rm Hom}_R (J,M/\sigma_\filtf(M)).\]
It is not obvious that $Q_\filtf(M)$ is naturally 
a left $R$-module: the fact that $\filtf$ is a Gabriel filter plays
a crucial role. Namely, given $f \in Q_\filtf(M)$, choose $J \in \filtf$
such that there is a $f_J$ in 
${\rm Hom}_R (J,M/\sigma_\filtf(M))$ representing $f$. 
For $r \in R$, the left ideal $(J:r)\in \filtf$, by the definition of a uniform
filter, and the rule $x \mapsto f(xr)$ defines an element $(rf)_{(J:r)}$  in
${\rm Hom}_R ((J:r),M/\sigma_\filtf(M))$ representing the class of
$rf$. This yields a well defined left action. 

\vskip .003in
Left multiplication by an element $r \in R$ defines a class
$[r] \in Q_\filtf(R)$. There is a unique ring
structure on $Q_\filtf(R)$, such that the
correspondence $i_\filtf : r \mapsto [r]$
becomes a ring homomorphism $i_\filtf : R \rightarrow Q_\filtf(R)$.

\vskip .011in
Not only every Gabriel filter defines an idempotent radical,
every radical also defines a Gabriel filter by the rule
\[  \filtf_\sigma = \{ J \subset R\,|\, \sigma(J) = 0\}.\]
When we restrict to the idempotent radicals, then this rule
gives a bijection between the idempotent radicals
and Gabriel filters.

\vskip .013in
Though it does not behave as nicely as Ore localization does,
scarcity of Ore sets makes Gabriel localization 
attractive and it is widely used. 
Moreover, this more general class of localizations
can be phrased fully in the language
of the Abelian category of left $R$-modules, and it generalizes to
other Abelian categories with some good properties. 
A common generality in which this is studied are
Grothendieck categories.
A {\bf Grothendieck category} ${\cal A}$ is an Abelian category
which is {\it cocomplete} (small inductive limits always exist
~\cite{Borceux,MacLane,vOyst:refl}),
where filtered limits are exact, and which 
posses a {\bf generator} (object $G$ in ${\cal A}$
such that $C \mapsto {\rm Hom}_{\cal C}(C,G)$ is a faithful functor),
for details 
cf.~\cite{JaraVerschorenVidal,Popescu,stenstrom,vOyst:refl, Smith:nag}.
Such categories are a natural
place to study noncommutative algebraic geometry beyond
the affine and 
projective cases~\cite{Rosen:book,Smith:nag,NdirvOyst:GroRep}.

\vskip .005in
A {\bf thick subcategory} of an Abelian category ${\cal A}$
is a replete (= full and closed under isomorphisms) 
subcategory ${\cal T}$ of ${\cal A}$ which is closed under
extensions, subobjects and quotients. In other words, 
an object $M'$ in a short sequence 
$0\rightarrow M \rightarrow M' \rightarrow M''\rightarrow 0$
in ${\cal A}$ belongs to ${\cal T}$ iff $M$ and $M''$ do. 
Localization at thick subcategories is a common framework in 
noncommutative algebraic geometry~\cite{Rosen:book,Ros:localalg}.
Starting from a pair $({\cal A}, {\cal T})$ where ${\cal A}$ is
Abelian and ${\cal T}$ is thick, one forms
a (Serre) quotient category~\cite{Borceux,Gab:catab,GZ,JaraVerschorenVidal,
Popescu,Ros:localalg}. As objects one takes
the objects of the original category, but in addition to the 
original morphisms one adds to the class of morphisms the formal inverses
of those morphisms $f$ for which both
${\rm Ker}\,f$ and ${\rm Coker}\,f$ are in ${\cal T}$.
A thick subcategory is called (Serre) {\bf localizing subcategory}
if the morphisms which are {\it invertible} in the quotient
category are exactly those for which 
${\rm Ker}\,f$ and ${\rm Coker}\,f$ are both in ${\cal T}$.
Hence, more than one thick subcategory may give the same
quotient category, and that ambiguity is removed if
we consider the corresponding localizing subcategories instead.

\vskip .01in
For any idempotent radical $\sigma$ in ${\cal A}$, 
define the class ${\cal T}_\sigma$ of $\sigma$-{\bf torsion objects}
and the class ${\cal F}_\sigma$ of $\sigma$-{\bf torsion free objects} by
\[ {\cal T}_\sigma = \{ M \in {\rm Ob}\, {\cal A}\,|\, \sigma(M) = M\},
\,\,\,\,\,\,\,
{\cal F}_\sigma = \{ M \in {\rm Ob}\, {\cal A}\,|\, \sigma(M) = 0\}.\]
Pair $({\cal T}_\sigma,{\cal F}_\sigma)$
is an example of a torsion theory and ${\cal T}_\sigma$
is a thick subcategory of ${\cal A}$. A {\bf torsion theory}
~\cite{Borceux,JaraVerschorenVidal}
in Abelian category ${\cal A}$ is a pair $({\cal T}, {\cal F})$
of replete subcategories of ${\cal A}$ such that
every morphism $T \rightarrow F$, where $T$ is an object in ${\cal T}$
and $F$ an object in ${\cal F}$, is zero morphism; and such that
every object $A$ in ${\cal A}$ can be placed into exact
sequence $0 \rightarrow T \rightarrow A \rightarrow F\rightarrow 0$
where $T$ is an object in ${\cal T}$ and $F$ an object in ${\cal F}$.
Localization in Abelian categories, and categories of modules
in particular, is often conveniently described in terms of
torsion theories.

\vskip .01in
Torsion theory is hereditary if every submodule of a torsion
module is also torsion. Torsion theories correspond to
idempotent preradicals. Hereditary
torsion theories correspond to Gabriel localizations,
which in turn correspond to idempotent radicals.

\section*{5. Covers via localizations.} 
\cite{Rosen:88,Ros:NcSch,JaraVerschorenVidal, Skoda:locelec}
For a moment, we take the most general view~\cite{Ros:NcSch,GZ,Borceux} 
that a {\bf localization} is a functor $Q^* : {\cal C}\rightarrow {\cal C}'$,
which is universal~\cite{GZ} with respect to inverting 
some class $\Sigma$ of morphisms in ${\cal C}$.
A functor $Q^*$ between {\it Abelian categories} 
will be called {\bf continuous}
if it has a right adjoint, say $Q_*$, and {\bf flat} if 
$Q^*$ is, in addition, exact.
A characterization of a localization functor is~\cite{GZ}
that it has a fully faithful right adjoint.
A family ${\cal V} = 
(Q^*_\lambda : {\cal C} \rightarrow {\cal C}_\lambda)_{\lambda \in \Lambda}$
of flat functors {\bf covers} ${\cal C}$ if it is {\it conservative} i.e.
\[ \forall f \in {\rm Mor}\,{\cal C} \,\,\,
\left(\,(\,(\forall \lambda)\,\,\,\,Q_\lambda^*(f) \,\mbox{ is invertible})\,
\,\,\Rightarrow\,\,
f \,\mbox{ is invertible.}\,\right)\]
We are interested here in covers by localizations only,
but we expect more general flat covers to play role in
future extensions of this work, as they do in commutative
algebraic geometry. Flat covers by localizations 
{\sc A.~Rosenberg} calls~\cite{Ros:NcSch}
{\bf Zariski covers}.

\vskip .014in
Let us now specialize to the category of left $R$-modules. 
We would like to carry further the picture that localizations correspond to
certain open subsets. In addition to covers,
one would like to have ``intersections''. A newcomer to localization
should be warned, however, of pitfalls here.

\vskip .013in
In the case of left Ore localizations $S^{-1}R$ and $T^{-1}R$, 
the natural candidate for intersection is the localization
at the set $S \vee T$ multiplicatively generated by $S$ and $T$. 
It is automatically Ore, hence $(S \vee T)^{-1} R$ is a ring as usually.
The set $ST$ of products $st$, $s \in S, t \in T$ is, 
in general, not multiplicative, but it does satisfy the left Ore condition.

\vskip .01in
As $S^{-1}R$ is an $R$-module, one
can always introduce $T^{-1}S^{-1}R$ as an $R$-module, by
applying the localization at $S$ first,
and the localization at $T$ after that.
$T^{-1}S^{-1}R$ is not necesarily a ring via Ore construction, 
as inverting $T$ (more precisely, $i_S(T)$) in $S^{-1}R$ 
by Ore method asks for $i_S(T)$ to be left Ore 
in $S^{-1}R$ what is not true in general,
and replacing left Ore sets by 2-sided Ore sets does not help.
To get feeling for this phenomenon 
write down the Ore condition for $i_S(T)$ in $S^{-1}R$ and notice that
there is more to check than the original Ore conditions for $T$
in $R$ say. Similarily, we can consider
$S^{-1}T^{-1}R$ and even $T^{-1}S^{-1}T^{-1}R$ etc.

\vskip .008in
$(S \vee T)^{-1} R$ is isomorphic to $T^{-1}S^{-1}R$ as an $R$-module
precisely when $i_S(T)$ is left Ore in $S^{-1}R$, hence,
by symmetry, iff $i_T(S)$ is left Ore in $T^{-1}R$. If this is true,
what is rare in noncommutative case, we say that $S$ and $T$
are mutually compatible left Ore sets 
~\cite{JaraVerschorenVidal,vOyst:assalg}
(not to confuse with the compatibility
with coaction which is a central topic in this paper). 
Hence, if $S$ and $T$ are compatible, module $T^{-1}S^{-1}R$
has a natural ring structure (characterized also
by the requirement that $(i_{T},T^{-1}S^{-1}R)$
is an object in ${\cal C}(S^{-1}R, T)$). Following {\sc A.~Rosenberg},
we call a cover {\bf semiseparated} if the localizations
are pairwise compatible (this makes sense for more general
localizations than Ore).

\vskip .004in
As $S \vee T$ is a bigger set than $ST$, we loose some information 
(kill more ideals) using localization at $S \vee T$
instead of the consecutive localizations at $S$ and $T$.
More concretely, if we view $\vee$ as an operation
of taking ``intersection of open sets'' and if a ring
covered by localizations could be considered as ``union'' of such,
then we face the fact that ``intersection'' is not distributive
with respect to ``union'':
$\{(S \vee T_\lambda)^{-1} R\}_{\lambda \in \Lambda}$ does not
necesarily cover  $S^{-1} R$ if $\{T_\lambda^{-1}R \}_{\lambda \in \Lambda}$
covers $R$. 
\vskip .012in
However, there is a positive result which puts us in business:
\vskip .008in
{\bf Globalization lemma.} 
{\it Suppose a finite family of Ore localizations 
$\{T_\lambda^{-1}R \}_{\lambda \in \Lambda}$ covers $R$.
Then for every left $R$-module $M$ the sequence
\[  0 \rightarrow M \rightarrow \prod_{\lambda \in \Lambda} T_\lambda^{-1}M
\rightarrow \prod_{(\mu,\nu) \in \Lambda\times \Lambda} 
T_{\mu}^{-1} T_\nu^{-1} M \]
is exact, where the first morphism is $m \mapsto \prod i_\lambda(m)$
and the second is 
}
\[ \prod_\lambda m_\lambda \mapsto 
\prod_{(\mu,\nu)} (i^\mu_{\mu,\nu} (m_\mu) - i^\nu_{\mu,\nu} (m_\nu)). \]
Here $i_\lambda$ is the localization map and $i^\mu_{\mu,\nu}$ is
the natural map from $T_\mu^{-1} M$ to $T_\mu^{-1} T_\nu^{-1} M$,
hence the lower indices for $i$ denote the target and the order matters.

\vskip .006in
This statement has been generalizated for
Gabriel filters~\cite{Rosen:88}, cf. also~\cite{JaraVerschorenVidal}. 
This may be derived from the application of Barr-Beck's theorem in this
setup.

\vskip .013in
The meaning of the globalization lemma is that 
every module can be reconstructed by gluing from its localizations, 
provided the overlaps in successive localizations in both orders
are taken simultaneously into account, as it is in the general
picture of flat descent, and
this principle extends to triple etc. localizations. 
Two systematic methods to use this
basic fact about covers have been proposed. 

\vskip .014in
$1^{\rm st}$ method, proposed by {\sc F.~van Oystaeyen} 
and {\sc L.~Willaert}, is to organize covers of some considered type 
into a {\bf noncommutative} (analogue of a) {\bf Grothendieck
topology}~\cite{vOystWill:Groth, vOystWill:cohschematic,
JaraVerschorenVidal, NdirvOyst:GroRep}, noncommutativity refering to the fact
discussed above that the order in taking successive localizations
matters. The notion of a sheaf and
a quasicoherent sheaf are then directly defined 
in analogy to the commutative situation~\cite{vOyst:assalg}. 
Most of the work in this direction is focused on the case
of ${\Bbb Z}_{\geq 0}$-graded Noetherian algebras for which
a nontrivial finite cover by nontrivial Ore sets exist,
so called {\bf schematic algebras}~\cite{vOystWill:Groth,vOyst:assalg}.

\vskip .011in
$2^{\rm nd}$ method centers on  
a comonad~\cite{Ros:NcSch} associated to given flat cover,
to place it transparently into the general picture 
of flat descent~\cite{Beck, Del:tan, SGA1, Ros:NcSS}. Then one
associated a cosimplicial object~\cite{Ros:NcSch} to the comonad.
When applying various functors
to this construction the exactness properties of the functors
and of the comonad play the decisive role; 
the description of objects obtained by gluing local data 
depends on the applicability of Barr-Beck's 
theorem~\cite{Beck, Borceux, BarrWells, 
KontsRos, MacLane, MLMoerd:Sheaves}. 
In particular, this is suitable for descent-type questions, 
and the construction of quotients can also be understood that way.

\vskip .023in

We say that a family $T^{-1}_\mu R$ 
of Ore localizations {\bf naively covers} a ring $R$ if
\[  0 \rightarrow R \rightarrow 
\prod_{\lambda \in \Lambda} T_\lambda^{-1}R \rightarrow \prod_{\mu < \nu} 
(T_\mu \cup T_\nu)^{-1} R \]
is exact. Here the second map is \[ \prod_\lambda m_\lambda \mapsto 
\prod_{\mu < \nu} (i_{\mu,\nu} (m_\mu) - i_{\nu,\mu} (m_\nu))\]
where $i_{\mu,\nu}$ (omitted upper indices!) 
is the map 
$i_{\mu,\nu}: T_\mu^{-1}R \rightarrow (T_\mu \cup T_\nu)^{-1} R$.
Every naive cover is a cover, and every semiseparated cover is
a naive cover. For the case of two localizations only, 
naive covers coincide with semiseparated covers.
Covers appear more naturally than naive covers do,
and having a naive (but nonsemiseparated) cover 
does not guarantee much more than a cover can do.

\section*{6. Trivial principal bundles.}
Given an affine algebraic group $B$ 
with a regular right action $\nu: E\times B\rightarrow E$ 
on an affine variety $E$, define linear map
\[ \rho\equiv\rho_\nu\, : {\cal O}(E) \rightarrow {\cal O}(E \times B)
\cong {\cal O}(E)\otimes {\cal O}(B)\,\mbox{ by }\,
\rho_\nu(f)(e,b) = f(\nu(e,b)).\] 
Map $\rho_\nu$ is a right ${\cal O}(B)$-coaction and an algebra map. 
In noncommutative case ${\cal O}(B)$ will be generalized to arbitrary
Hopf algebra ${\cal B}$, and $E$ to any ${\cal O}(B)$-comodule
algebra ${\cal E}$.
Any left (right) comodule ${\cal E}$
over a bialgebra ${\cal B}$ such that ${\cal E}$ is an associative algebra
and the coaction is a homomorphism of algebras is called
a left (right) ${\cal B}$-{\bf {comodule algebra}}.

\vskip .01in
For a set $E$ with right $B$-action define subset
\[ E \times_B E = \{ (e_1,e_2) \in E \times E, \,|\,
\exists b \in B,\, e_1 b = e_2 \} \subset E \times E.\]
If $E$ is a topological space then $E\times_B E$ inherits
a subspace topology from $E \times E$. The action of $B$ is
{\bf free} if for every pair $(e_1,e_2) \in E \times_B E$
there is a {\it unique} $b$ with $e_1 b = e_2$. Then the rule
$\tau : (e_1,e_2) \mapsto b$ defines a map of sets
$\tau : E \times_B E \rightarrow E$. 
If $B$ is a {\it topological} group
a {\bf {principal $B$-bundle}} is a topological space $E$ 
with a free right $B$-action, such that $\tau$ is a continuous map.
In addition, a local triviality condition is usually required.

For algebraic varieties continuous maps are replaced by regular maps.

\vskip .01in
A principal $B$-bundle is {\bf trivial} if there is a
section $t : X \rightarrow E$ of the projection $p : E \rightarrow X$, 
i.e. a continuous map such that $p \circ t = {\rm id}_{X}$.
Let $t$ be a section, $p^{-1}(x)$ some fiber,
and $f$ a continuous function on $B$. Then the formula
\begin{equation}\label{eq:classical-gamma}
\gamma_t (f)(e) = f(\tau(t(p(e)),e))
\end{equation}
defines a continuous function $\gamma_t (f)$ on $E$.
In the algebraic case, map
\begin{equation}\label{eq:gamma-map}
 \gamma_t : {\cal B} \rightarrow {\cal E}
\,\,\,\,\,\,f \mapsto \gamma_t(f) 
\end{equation}
 defines
a map of commutative ${\cal B}$-comodule algebras, 
where ${\cal B} = {\cal O}(B)$ etc. 

\vskip .02in
To prepare for the study of locally trivial principal bundles
we now introduce certain free and smash products, and
a notion of compatibility.

\vskip .01in
Let ${\cal B}$ be a $\genfd$-bialgebra and 
$(M_\alpha,\rho_\alpha)$ a family of ${\cal B}$-comodules. 
A family of $\genfd$-linear maps $f_\alpha : M_\alpha \rightarrow A$
where $A$ is a fixed algebra is called 
{\bf {$\{\rho_\alpha\}$-compatible}}
iff there is a {\it unique} coaction $\rho_A$ on $A$
making $A$ a ${\cal B}$-comodule algebra and for each $\alpha$
the diagram
\[\begin{array}{ccl}
  M_\alpha &\stackrel{\rho_\alpha}{\rightarrow}& M_\alpha \otimes {\cal B}\\
  f_\alpha \downarrow & & \,\,\,\,\,\,\,\downarrow f_\alpha \otimes {\rm id} \\
  A  & \stackrel{\rho_A}{\rightarrow} & A \otimes {\cal B}
\end{array}
\]
commutes.
\remar{Examples.}{ 1. Compatible localization maps, cf.~8.1.

2. Given a right  ${\cal B}$-comodule $(V,\rho)$,
the inclusion $V \hookrightarrow T(V)$ is a $\rho$-compatible.

3. Let $(A_1, \rho_1), (A_2,\rho_2), \ldots, (A_n, \rho_n)$ 
be a finite sequence of right ${\cal B}$-comodule algebras. 
The family of natural inclusions $i_j : A_j \rightarrow
A_1 \otimes A_2 \otimes \ldots \otimes A_n$ is $\{\rho_j\}$-compatible.

4.  A free product of a family of associative algebras
$\{A_{j}\}_{ j \in {\cal J}}$
is $T( \oplus_{j \in \cal J} A_j )$ modulo an ideal
generated by all expressions of the form
$a \otimes a' - aa'$ where both $a,a'$ belong to $A_j$ with the same $j$
If all $A_j$ are ${\cal B}$-comodule algebras,
then combining Examples 2 and 3 and this description
conclude that the family of natural inclusions
$i_j : A_j\hookrightarrow \star_{j \in J} A_j$
is $\{\rho_j\}$-compatible.

4a. In particular,
consider an algebra $U$ with the trivial ${\cal B}$-coaction
and ${\cal B}$ self-coacting by comultiplication. By compatibility,
$U \star {\cal B}$ becomes a ${\cal B}$-comodule algebra.
}

Let $U$ be an algebra and ${\cal B}$ a Hopf algebra.
Define the {\bf category ${\cal C} = {\cal C}(U,{\cal B})$
of higher smash products of $U$ and ${\cal B}$} as follows. 
An object in ${\cal C}$ is a triple $(A,\iota,\gamma)$
where $A$ is a ${\cal B}$-comodule algebra, 
$\iota : U \hookrightarrow A^{{\rm co}\,{\cal B}}$ is 1-1 algebra map,
and  $\gamma  : {\cal B} \rightarrow A$ is a ${\cal B}$-comodule algebra map,
such that $A$ is generated by $\iota(U)$ and $\gamma({\cal B})$.
A morphism in ${\cal C}$ from $(A,\iota,\gamma)$ into
$(A',\iota',\gamma')$ is a map $f : A \rightarrow A'$ of ${\cal B}$-comodule
algebras such that $\iota' = \iota \circ f$ and $\gamma' = \gamma \circ f$.

\vskip .002in
For every object in ${\cal C}$, one defines a map
\[
\rhd : {\cal B} \otimes A \rightarrow A,\,\,\,\,\,\,\,\,b \rhd a = 
\sum\gamma(b_{(1)})a\gamma(Sb_{(2)}).
\]
This is an algebra action making $A$ a ${\cal B}$-module algebra
(i.e. $b \rhd 1_A = 1_A$ and 
$b \rhd (aa') = \sum (b_{(1)}\rhd a)(b_{(2)}\rhd a')$),
and $A^{{\rm co}\,{\cal B}}$ a ${\cal B}$-submodule subalgebra.
Moreover, $A^{{\rm co}{\cal B}}$ is the smallest 
${\cal B}$-stable (i.e. $\rhd$-invariant) subalgebra of $A$
containing $U$. It follows that $A$ is isomorphic, 
as a comodule algebra, to the (ordinary)
{\bf smash product} $A^{{\rm co}{\cal B}} \sharp \,{\cal B}$
for that action on $A^{{\rm co}{\cal B}}$, i.e. the tensor product
$A^{{\rm co}{\cal B}} \otimes {\cal B}$ with product rule
$(v \otimes b)(v'\otimes b') = v (b_{(1)}\rhd v') \otimes b_{(2)}b'$
and ${\cal B}$-coaction ${\rm id} \otimes \Delta_{\cal B}$.

\vskip .01in
The free product $U \star {\cal B}$ with the induced ${\cal B}$-coaction
(cf. Example 4a, above) and the natural inclusions
$\iota_U : U \hookrightarrow U \star {\cal B}$
and $\gamma : {\cal B} \hookrightarrow U \star {\cal B}$,
form a universal object in ${\cal C}$. 
Any ordinary smash product $U \sharp {\cal B}$
in ${\cal C}(U,{\cal B})$ is {\bf locally terminal}
in the sense that any morphism in ${\cal C}$ with source
$U \sharp {\cal B}$ is an isomorphism in ${\cal C}$,
and that for any object $C$ there is at most one morphism
from $C$ to $U \sharp {\cal B}$.
All locally terminal objects in ${\cal C}$ are of that type.
Isomorphism classes in ${\cal C}$ of ordinary smash products 
(with the same $U$)
are distinguished by action $\rhd$. For every given
${\cal B}$-module algebra action on $U$, there
is a unique isomorphism class of locally terminal objects 
in ${\cal C}$ such that $\iota_U$ is a map of ${\cal B}$-modules.

\vskip .009in
Finally, ${\cal C}$ is an {\bf umbrella category} i.e.
it posseses a universal initial object and a class $\theta$ of
locally terminal objects, and for every object $C$ in ${\cal C}$
there is at least one $T \in \theta$ with a (unique)
morphism $C \rightarrow T$ in ${\cal C}$.

\vskip .011in
A {\bf trivial quantum principal ${\cal B}$-bundle} is
an object in ${\cal C}(U,{\cal B})$ for some algebra $U$.
Notice that {\it for the same ${\cal B}$ and the same
base space $U$ there is more than one trivial bundle},
as the action $\rhd$ can be different. 
In our point of view even for fixed total space we allow 
different $i(U)\subset A^{{\rm co}{\cal B}}$ 
as long as $\iota(U)$ and $\gamma({\cal B})$ generate $A$.
Namely a natural {\it candidate} for the quantum base space 
is sometimes smaller than the whole algebra of coinvariants,
and the latter appears with the help of action in total
space (so it is not fully ``base'').

\section*{7. Commutative local triviality and torsors.}
For algebraic principal bundles,
local triviality is considered with respect to one
of the several standard topologies for schemes. 
Local triviality in Zariski topology is the strongest
requirement, and the local triviality in 
\'{e}tale, {\it fppf}, {\it fpqc} 
topology are weaker, in that order~\cite{Dem-Gab-eng, Milne}. 
A principal bundle locally trivial in \'{e}tale topology is often called
a {\bf torsor}. If the orbit space is denoted by $X$, one
can replace $B$ by trivial $B$-bundle ${\bf B}$ over $X$
(i.e. by the product $B\times X$). Then  ${\bf B}$ is a (relative) group scheme
over $X$. More generally,  
consider any group scheme ${\bf B}$ over $X$ 
(topological analogue: a bundle of groups over $X$)
acting upon $X$ in category of schemes over $X$. 
Descent theory implies~\cite{Milne} that {\it local triviality of $E$ over $X$
in flat topology is equivalent to the requirement that $E$ is
faithfully flat and locally of finite-type over $X$ and} that 
\begin{equation}\label{eq:torsor-map}
 (e,b) \mapsto (e,eb) : E \times_X {\bf B} \rightarrow E \times_X E 
\end{equation}
{\it is an isomorphism.} If $E,X$ and $B$ are {\it affine}, then
we can dualize~(\ref{eq:torsor-map}) 
by taking global sections
of the structure sheaf~\cite{Dem-Gab-eng}. In that case
${\cal X} = \Gamma {\cal O}_X$ is the ring of $\Gamma {\cal O}_B$-coinvariants
in ${\cal E} = \Gamma {\cal O}_E$. 
Then, if~(\ref{eq:torsor-map}) is isomorphism, 
it follows~\cite{Dem-Gab-eng} that the map
\begin{equation}\label{eq:Hopf-Galois}
 {\cal E} \otimes_{\cal X} {\cal E} \rightarrow 
{\cal E} \otimes_\genfd {\cal B},\,\,\,\,\,\,\,
 e \otimes_{\cal X} e' \mapsto (e \otimes_\genfd 1_{\cal B}) \rho(e')
\end{equation}
is bijective. An {\bf extension} of any algebra ${\cal X}$
by a Hopf algebra ${\cal B}$ is a ${\cal B}$-comodule algebra ${\cal E}$
such that ${\cal X}$ equals to the ring ${\cal E}^{{\rm co}{\cal B}}$
of coinvariants in ${\cal E}$.
An extension is {\bf Hopf-Galois} iff the map~(\ref{eq:Hopf-Galois})
is bijective. 
An extension is {\bf cleft} if there is
a convolution-invertible map of ${\cal B}$-comodules
$\gamma : {\cal B} \rightarrow {\cal E}$
(cf.~(\ref{eq:gamma-map})). For any cleft extension
there is a $\genfd$-linear map
\[
\rhd : {\cal B} \otimes {\cal X}
\rightarrow {\cal X},\,\,\,\,\,\,\,\,b \rhd u = 
\sum\gamma(b_{(1)})u\gamma^{-1}(b_{(2)}),
\]
as it is direct to check that the right-hand side 
is a ${\cal B}$-coinvariant. Map $\rhd$ {\bf measures} ${\cal X}$
i.e. $b \rhd 1 = 1$ and $b \rhd (uv) = \sum (b_{(1)}\rhd u)(b_{(2)}\rhd v)$. 
Cleft extensions are (a special case of) Hopf-Galois extension.
$({\cal E},\id_{{\cal E}^{{\rm co}\cal B}},\gamma)$
is a {\bf trivial quantum principal ${\cal B}$-bundle} 
over ${\cal E}^{{\rm co}B}$ if $\gamma$ is also an algebra map. 
In that case, convolution invertibility
of $\gamma$ comes for free as $\gamma^{-1} = \gamma \circ S$. 

\section*{8. Localized coinvariants.}
{\bf 8.1} {\it Compatibility.}
Suppose we are given a bialgebra $\cal{B}$ and
a right $\cal{B}$-comodule algebra $\cal{E}$.
An Ore set $T$ and localization
$i_T : {\cal E} \rightarrow T^{-1}{\cal E}$
are {\bf compatible} with coaction
$\rho : {\cal E} \rightarrow {\cal E}\otimes {\cal B}$
if there exist a unique map
\[ \rho_T : T^{-1}{\cal E} \rightarrow 
T^{-1}{\cal E}\otimes {\cal B} \]
which makes $(T^{-1}{\cal E},\rho_T)$ into
a $\cal{B}$-comodule algebra and the following diagram commutes
\[\begin{array}{ccl}
{\cal E}&\stackrel{\rho}{\rightarrow} & 
\,\,\,\,\,\,\,\,\,\,\,{\cal E}\otimes {\cal B}\\
i_T \downarrow && i_T \otimes {\rm id} \downarrow\\
T^{-1}{\cal E}&\stackrel{\rho_T}{\rightarrow}\,
&\,\,\,\,\,\,T^{-1}{\cal E}\otimes {\cal B}
\end{array}\]
i.e. $\rho_T \circ i_T = (i_T \otimes {\rm id})\circ \rho$.

In other words, Ore localization is $\rho$-compatible iff there is a unique
$\cal{B}$-comodule algebra structure on the localized algebra
such that the localization map is an intertwiner.

\vskip .01in
This definition is still appropriate, for 
those more general (than Ore) localizations,
for which $Q(R)$ is still naturally a ring
with homomorphism $i : R \rightarrow Q(R)$, {\it and}
the composition of functors $Q^* i_*$ (where
$i_*$ is the restriction of scalars functor
from $Q(R)$--modules to $R$--modules)
is an equivalence between the category of
$Q(R)$-modules and the localized category. (Only) such
localization functors satisfy $Q(M) = Q(R) \otimes_R M$ for all $M$. 
They are called {\bf perfect localizations}
~\cite{JaraVerschorenVidal,stenstrom}. A Gabriel localization is
perfect iff $Q = Q_* Q^*$ is an {\it exact} endofunctor.
For even more general cases, one may 
redefine the $\rho$-compatibility in the language of Hopf module
categories.

\vskip .01in
Classically we think of $\rho$-compatibility
as the condition that the corresponding
Zariski open set is $B$-invariant, i.e. a union of $B$-orbits.

\vskip .01in
Any $\rho$-compatible Ore localization
$i_S : {\cal E} \rightarrow S^{-1}{\cal E}$ of a cleft ${\cal B}$-extension
${\cal E}$ is cleft.  The section is $\gamma_S = i_S \circ\gamma$.
If $\gamma$ is an algebra map, so is $\gamma_S$. This is forced by
the very definition of the $\rho$-compatibility.

\vskip .02in
{\bf 8.2}
{\it Practical criterium of compatibility.} Localization ${\cal E}[T^{-1}]$
is $\rho_T$-compatible iff for each $t \in T$, element
$((i_T \otimes {\rm id}) \circ \rho)(t)$ is invertible
in algebra $T^{-1}{\cal E} \otimes {\cal B}$.

The proof is elementary, cf.~\cite{Skoda:thes,Skoda:loc-coinv1}.

\vskip .02in
{\bf 8.3} {\it Localized coinvariants.} ~\cite{Skoda:thes,Skoda:loc-coinv1}
Let $T$ be a $\rho$-compatible right Ore set
in a right  ${\cal B}$-comodule algebra $({\cal E}, \rho)$.
By compatibility there is a
uniquely defined localized ${\cal B}$-coaction $\rho_T$
on ${\cal E}_T = T^{-1}{\cal E}$. 
We define the {\bf algebra of} $T$-{\bf localized} 
right $\rho$-{\bf coinvariants} in ${\cal E}$
to be the algebra of $\rho_T$-coinvariants in ${\cal E}_T$:
 \[ {\cal E}_T^{{\rm co}\cal B} = \{ y \in {\cal E}_T \,|\,
\rho_{T} y = y \otimes 1 \}. \]

\vskip .02in
{\bf 8.4} {\it Nested localizations.}~\cite{Skoda:thes,Skoda:loc-coinv1}
Let $S \subset T \subset {\cal E}$ be an inclusion of
Ore subsets in ${\cal E}$, both compatible with ${\cal B}$-coaction $\rho$.

(i) The square diagram involving localized coactions commutes
\[\begin{array}{ccl}
S^{-1}{\cal E} & \stackrel{\rho_S}{\rightarrow} & 
        S^{-1}{\cal E}  \otimes\, {\cal B}\\
 i^S_{T} \downarrow & & \,\,\,\,\,\,\,\,\,\,\,\,\,\,\,
  \downarrow i^S_{T} \otimes {\rm id}\\
T^{-1}{\cal E} & \stackrel{\rho_T}{\rightarrow} & 
        T^{-1}{\cal E}  \otimes\, {\cal B}
\end{array}\]
In other words, the natural maps $i^S_T$ 
between localizations are intertwiners.

(ii) The natural map $i^S_{T}$ between the localizations maps 
the (sub)algebra of $S$-localized coinvariants into 
the (sub)algebra of $T$-localized coinvariants:
\[ i^S_{T}(({S^{-1}{\cal E}})^{{\rm co}\cal{B}}) \subset 
({T^{-1}{\cal E}})^{{\rm co}{\cal B}}.\]

\vskip .02in
{\bf 8.5} {\it Compatibility and Hopf modules.}
Given a right ${\cal B}$-comodule algebra ${\cal E}$, 
(relative) $({\cal E}, {\cal B})$--{\bf Hopf modules} 
or simply dimodules, are the objects of the
category ${}_{\cal E}{\cal M}^{\cal B}$
formed by left ${\cal E}$-modules $M$, which are also right 
${\cal B}$-comodules, with compatibility 
\[ \sum e_{(0)} m_{(0)}\otimes_\genfd e_{(1)} m_{(1)} = 
        \sum (em)_{(0)} \otimes_\genfd (em)_{(1)}.\]
For any $M \in {}_{\cal E}{\cal M}^{\cal B}$ one defines
the ${\cal B}$-coaction $\rho'$ 
on the tensor product $S^{-1}{\cal E} \otimes_\genfd M$
over ground field $\genfd$:
\[\begin{array}{l}
\rho' : S^{-1}{\cal E} \otimes_\genfd M \rightarrow 
S^{-1}{\cal E} \otimes_\genfd M \otimes_\genfd {\cal B}
\\
\rho'(x \otimes_\genfd m) 
= \sum (x_{(0)} \otimes_{\genfd} m_{(0)}) 
\otimes_\genfd x_{(1)} m_{(1)}.
\end{array}\]
This coaction factors to the $\cal B$-coaction $\rho_M$ on the
$S^{-1}M$ which is the tensor product of the same factors 
but now over ${\cal E}$:
\[
  \rho_{S^{-1}M}(x \otimes_{\cal E} m) = 
\sum (x_{(0)} \otimes_{\cal E} m_{(0)}) 
\otimes_\genfd x_{(1)} m_{(1)}
\]
i.e. we have the commutative diagram
\[\begin{array}{ccc}
\,\,\,\,S^{-1}{\cal E} \otimes_\genfd M & \stackrel{\rho'}{\rightarrow} &
(S^{-1}{\cal E} \otimes_\genfd M) \otimes_\genfd {\cal B} \\
\downarrow {\rm pr} & & \downarrow {\rm pr} \otimes_\genfd {\rm id}   \\
S^{-1}M = S^{-1}{\cal E} \otimes_{\cal E} M 
& \stackrel{\rho_{S^{-1}M}}{\rightarrow} &
(S^{-1}{\cal E} \otimes_{\cal E} M) \otimes_\genfd {\cal B},
\end{array}\]
where the vertical maps are the natural projections. The bottom
map is well-defined
thanks to the compatibility of ${\cal E}$-module and
${\cal B}$-comodule structure on $M$:
\[\begin{array}{l} \sum ((ye)_{(0)} \otimes_{\cal E} n_{(0)}) \otimes_\genfd 
(ye)_{(1)} n_{(1)} =
\\ \,\,\,\,\,\,\,\,\,\,\,\,\,\,\,\,\, 
=\sum (y_{(0)} e_{(0)} \otimes_{\cal E} n_{(0)}) \otimes_\genfd
y_{(1)} e_{(1)} n_{(1)}
\\ \,\,\,\,\,\,\,\,\,\,\,\,\,\,\,\,\, 
=\sum (y_{(0)}\otimes_{\cal E} (en)_{(0)}) \otimes_\genfd
y_{(1)} (en)_{(1)} \in S^{-1}M \otimes_\genfd {\cal B}.
\end{array}\]

Thus we get a functor $Q^{{\cal B}}_S : {}_{\cal E}{\cal M}^{\cal B}
\rightarrow {}_{S^{-1}{\cal E}}{\cal M}^{{\cal B}}$.

This argument extends to perfect localizations.
For more general localizations, when the
definition in terms of the ground ring is not appropriate,
we say that a localization functor $Q = Q^*$ is $\rho$-compatible 
if it induces a localization functor $Q^{\cal B}$
on the category of dimodules ${}_{\cal E}{\cal M}^{{\rm co}\cal B}$. 
The localization maps $i_M : M \rightarrow Q(M)$ 
then also lift to the maps of dimodules.
This may be seen by observing that, abstractly,
the localization maps come from 
the adjunction $i : {\rm Id} \rightarrow Q_* Q^*$.

Apart from easy generalizations, the reformulation in terms of
dimodule categories has other useful consequences.
For each dimodule $M$, there is an equality
$Q^{\cal B}_S(Q^{\cal B}_S M) = Q^{\cal B}_S M$ of dimodules.
More important observation is
that the functors of the type $Q^{\cal B}_S$ can be iterated.
This means that there is a natural ${\cal B}$-coaction 
on the successive $\rho$-compatible localizations 
$S^{-1}_1 S^{-1}_2 \cdots S^{-1}_n {\cal E}$
what is a refinement (cf. discussion of covers
by localizations above) of the previous picture where
we could do this only for  
$(S_1 \vee S_2 \vee \ldots \vee S_n)^{-1} {\cal E}$,
as the latter is a ring. We can now define the {\bf module
of localized coinvariants} in $S^{-1}_1 S^{-1}_2 \cdots S^{-1}_n {\cal E}$
as the module of coinvariants for this induced ${\cal B}$-coaction. 

\section*{\bf 9. Covers by coaction--compatible localizations.}
To every cover of ${\cal E}$ by flat localizations 
one associates a cosimplicial object
\begin{equation}\label{eq:cover-res}
 \diagram{  
{\cal E} \fhd{}{} \prod_\mu {\cal E}_\mu \fhd{\fhd{}{}}{} \prod_{\mu,\nu}
{\cal E}_{\mu\nu} \fhd{\fhd{}{}}{\fhd{}{}} \,\,\,\cdots
}
\end{equation}
where ${\cal E}_{\mu\nu...}$ are the successive localizations.
If the cover is semiseparated then we can identify ${\cal E}_{\mu\nu}$
with  ${\cal E}_{\nu\mu}$ etc. so the products
can be taken with $\mu < \nu < \ldots$ with respect to any fixed
ordering on the indices of the cover. If the
localizations are $\rho$-compatible and cover semiseparated 
thia is a cosimplicial object in the category of ${\cal B}$-comodule
algebras. Without semiseparetedness
it is only in category ${}_{\cal E}{\cal M}^{\cal B}$.

\vskip .02in
Denote by ${}_{\cal E}{\cal M}_{\mu\nu\ldots}^{\cal B}$ the
category obtained from ${}_{\cal E}{\cal M}^{\cal B}$ by
successive application of dimodule category
localizations $Q^{\cal B}_\mu$. Again, in semiseparated case,
${}_{\cal E}{\cal M}_{\mu\nu\ldots}^{\cal B}$  will agree
with the dimodule category ${}_{{\cal E}_{\mu\nu\ldots}}{\cal M}^{\cal B}$,
but in general ${\cal E}_{\mu\nu\ldots}$ is not a ring, but
only an ${\cal E}$-bimodule, so the latter category does not make
sense. For any flat localization functor $Q^*$, 
its right adjoint $Q_*$ is fully faithful ~\cite{GZ}. 
Thus we may view the objects in
${}_{\cal E}{\cal M}_{\mu\nu\ldots}^{\cal B}$ as living
in ${}_{\cal E}{\cal M}^{\cal B}$. In particular, every functor
$F$ defined on ${}_{\cal E}{\cal M}^{\cal B}$ can be
naturally evaluated on ${}_{\cal E}{\cal M}_{\mu\nu\ldots}^{\cal B}$.

\vskip .012in
More generally, if we replace 
${\cal B}$-comodule algebra ${\cal E}$ in~(\ref{eq:cover-res})
by any dimodule $M \in {}_{\cal E}{\cal M}^{\cal B}$
we also obtain a cosimplicial object. Moreover,
the construction is functorial in $M$. 
Hence we obtain a cosimplicial object 
in the category of endofunctors 
${\rm End}\,{}_{\cal E}{\cal M}^{\cal B}$.

\vskip .007in
For details on this ``standard resolution'' 
construction see~\cite{Ros:NcSch},
remembering that now our objects live in Abelian (actually
Grothendieck~\cite{wisbSurvey}) category ${}_{\cal E}{\cal M}^{\cal B}$.

\vskip .02in
Consider the left exact functor 
\[ ()^{{\rm co}{\cal B}} : {}_{\cal E}{\cal M}^{\cal B}
\rightarrow {}_{\cal U}{\cal M}\]
of taking submodule of coinvariants, where $U = {\cal E}^{{\rm co}{\cal B}}$.
Functor $M \mapsto {\cal E}\otimes_{\cal U} M$ is a left adjoint
to $()^{{\rm co}{\cal B}}$. Here the ${\cal E}$-action on  
${\cal E}\otimes_{\cal U} M$ extends map $e(e' \otimes_{\cal U} m) 
= ee' \otimes_{\cal U} m$
and ${\cal B}$-coaction extends map
$e \otimes m \mapsto \sum
(e_{(0)} \otimes_{\cal U} v) \otimes_\genfd e_{(1)}$.
The adjunctions are given by
\[\begin{array}{ll} 
M \mapsto ({\cal E} \otimes_{\cal U} M)^{{\rm co}{\cal B}},\,
& m\mapsto 1 \otimes m,\\
{\cal E} \otimes_{\cal U} M^{{\rm co}{\cal B}} \rightarrow M,\,
& e \otimes_{\cal U} m \mapsto em.
\end{array}\]
{\sc H.-J. Schneider}{\bf 's theorem}~\cite{Schn-prinHomSp,Nuss:descent}
says that {\it these two adjoint functors are equivalences iff
${\cal E}$ is a faithfully flat Hopf Galois extension of ${\cal U}$.}
We'll now sketch interplay between the functors playing role
in the Schneider's theorem and the (co)simplicial structures
coming from the cover by the coaction-compatible localizations.
\vskip .02in
Throughout this section we denote ${\cal U}^{\mu\nu\ldots} =
({\cal E}_{\mu\nu\ldots})^{{\rm co}{\cal B}}$. We use the upper
indices, as  ${\cal U}^{\mu\nu\ldots}$ is not necessarily
a localization (lower indices!) of ${\cal U}$.
In section 8, we have discussed that any $\rho$-compatible localization
functor $Q_\mu$ induces the localization functor
$Q^{\cal B}_{\mu}$ on the category of dimodules,
hence if the functors from the Schneider's theorem are equivalences
then there is a well defined
localized category ${}_{\cal U}{\cal M}_{\mu\nu\ldots}$
of ${}_{\cal U}{\cal M}$, but in general ${}_{\cal U}{\cal M}$ is not defined. 
When such localized categories are defined,
we can restate the very fact by saying that
the $\rho$-compatible cover 
${\frak U} = \{ Q_\mu \}$ of ${\cal E}$ induces
a cover $\tilde{\frak U}$ of ${\cal U}$.
The category of
${\cal U}^{\mu\nu\ldots}$-modules, denoted by 
${}_{{\cal U}^{\mu\nu\ldots}}{\cal M}$ is not necessarily
a localized category ${}_{\cal U}{\cal M}_{\mu\nu\ldots}$
of ${}_{\cal U}{\cal M}$, even when the latter
exists (when it exists it is typically bigger than the latter). 
In other words, taking the coinvariants and localizing
do not commute.
\vskip .02in
A cosimplicial object in ${}_{\cal E}{\cal M}^{\cal B}$
is given by a functor from the cosimplicial
category $\Delta$ to ${}_{\cal E}{\cal M}^{\cal B}$.
In our case, to each function $f : [l] \rightarrow [k]$ 
in the cosimplicial category $\Delta$ 
one assigns a certain composition $Q^{\cal B}_f$ of
the localization functors $Q_\mu^{\cal B}$ 
(roughly speaking: formulas also involve their right adjoints,
the inclusion into the direct product, and the unit of the
appropriate adjunction, cf.~\cite{MacLane,Ros:NcSch,Ros:NcSS,Weibel}).

\vskip .03in
For each cover ${\frak U} = \{ Q_\mu \}$ of ${\cal E}$ by
$\rho$-compatible Gabriel localizations we define category
${}_{\cal E}{\cal M}({\frak U})^{{\rm co}\cal B}$, 
the ${\frak U}$-{\bf accessible quotient}, as follows. 
For mere sake of comparison, we'll also define two other
categories in similar way, ${}_{\cal E}{\cal M}^{\cal B}({\frak U})$,
${}_{\cal U}{\cal M}(\tilde{\frak U})$.
Instead of an object $M$ 
and its successive localizations
$M_{\mu\nu\ldots}$ we consider
a family $N^\cdot$ of objects $N^{\mu\nu\ldots}$
({\it upper} indices!)
in categories ${\cal C}^{\mu\nu\ldots}$ respectively,
given with a {\it family of structure morphisms} $b_\cdot$.  
Here ${\cal C}^{\mu\nu\ldots}= {}_{{\cal U}^{\mu\nu\ldots}}{\cal M}$
for the ${\frak U}$-accessible quotient 
${}_{\cal E}{\cal M}({\frak U})^{{\rm co}\cal B}$, and 
${\cal C}^{\mu\nu\ldots} ={}_{\cal E}{\cal M}^{\cal B}_{\mu\nu\ldots}$
and ${}_{\cal U}{\cal M}_{\mu\nu\ldots}$ for the construction of
${}_{\cal E}{\cal M}^{\cal B}({\frak U})$ and
${}_{\cal U}{\cal M}(\tilde{\frak U})$ respectively.
Here ${}_{\cal U}{\cal M}(\tilde{\frak U})$ is defined only
when $\tilde{\frak U}$ is defined.
Denote $N_k = \prod_{\mu_1,\ldots,\mu_k} N^{\mu_1\ldots\mu_k}$
and ${\cal C}_k$ the corresponding category.
A {\bf family of structure morphisms} is 
an assignment $b: (f,k) \mapsto b(f,k)$,
where $f : [l] \rightarrow [k]$ is in the cosimplicial
category $\Delta$, $k$ is a natural number and
$b(f,k)$ is a morphism in ${\cal C}_k$,
\[     b(f,k) : Q_f^{\cal B}(N_k) \rightarrow  N_{l},
\mbox{ such that }\left\lbrace\begin{array}{l}
 b(f \circ g,k) = b(g,l)\circ Q_g^{\cal B}(b(f, k))
\\
Q_{f \circ g}^{\cal B}(N_k) =  
Q_g^{\cal B} Q_f^{\cal B}(N_k) \rightarrow 
Q_g^{\cal B} (N_l) \rightarrow N_m\\
b({\rm id},k) = {\rm id}
\end{array}\right.
\] 
Objects of the category ${}_{\cal E}{\cal M}({\frak U})^{{\rm co}\cal B}$ 
are the pairs $(N_\bullet, b)$, as described above, and similarily for
the other two categories. A morphism
$f_\bullet :(N_\bullet,b)\rightarrow (N'_\bullet,b')$
is a sequence $(f_k)$, where each $f_k : N_k \rightarrow N'_k$ is
an ${\cal E}_k$-module map, and were $f_k$ commute with structure
morphisms. Composition of morphisms is defined componentwise.
This construction
is an analogue of a
category of simplicial sheaves over a simplicial space 
cf.~\cite{BernLunts,Del:HodgeIII,Jard:simp,Moerd:class}.
\vskip .01in
Functor $M \mapsto (N_\bullet, b)$ 
defined by $N^{\mu\nu\ldots} = (M_{\mu\nu\ldots})^{{\rm co}{\cal B}}$
and $b = {\rm id}$ defines a functor
${}_{\cal E}{\cal M}^{\cal B}\rightarrow
{}_{\cal E}{\cal M}({\frak U})^{{\rm co}\cal B}$.
Functor $M \mapsto (N_\bullet, b)$ 
defined by $N^{\mu\nu\ldots} = M_{\mu\nu\ldots}$
from ${}_{\cal E}{\cal M}^{\cal B}$ to 
${}_{\cal E}{\cal M}^{\cal B}({\frak U})$
is, however, an equivalence of categories 
(and all the structure morphisms are isomorphisms). 
This follows from the
globalization lemma (section 5) and general nonsense about
simplical objects (cf.~\cite{Ros:NcSS}, Prop. 1.0.8.2, for
the appropriate statement in the language of comonads).
\vskip .01in
It is a basic phenomenon that taking coinvariants does not
commute with localization. However, if $()^{{\rm co}{\cal B}}$
is an equivalence of categories, then all the localization
functors $Q^{\cal B}_\mu$ may be viewed as the localization
functors in category ${}_{\cal U}{\cal M}$. By the
globalization lemma (now in ${}_{\cal U}{\cal M}$), 
${}_{\cal U}{\cal M} = {}_{\cal U}{\cal M}(\tilde{\frak U})$,
where $\tilde{\frak U}$ is the localization cover induced by
${\frak U}$ as above. Hence ${}_{\cal E}{\cal M}^{\cal B}$
is equivalent to category ${}_{\cal E}{\cal M}({\frak U})^{{\rm co}\cal B}$.
If $()^{{\rm co}{\cal B}}$ is not an equivalence, 
but  $()^{{\rm co}{\cal B}}$ of each of the
categories ${}_{\cal E}{\cal M}^{\cal B}_\mu$ is
(each ${\cal E}_\mu$ is a faithfully flat Hopf-Galois extension),
then ${}_{\cal E}{\cal M}^{\cal B}$
is still equivalent to ${}_{\cal E}{\cal M}({\frak U})^{{\rm co}\cal B}$.

\vskip .02in
One thinks of the category of Hopf modules ${}_{\cal E}{\cal M}^{\cal B}$
as the category $Qcoh_B(E)$ of $B$-equivariant 
(for right $B$-action, 
where $B$ is an affine group and ${\cal B} = {\cal O}(B)$) 
quasicoherent sheaves on affine $B$-variety $E$.
If the action is free and under some flatness conditions, classical
descent theory~\cite{SGA1}
identifies $Qcoh_B(E)$ with category $Qcoh(E/B)$ 
(in that case, 'descent is effective')
of quasicoherent sheaves on the quotient space.
The latter category is approximated by the ${\frak U}$-accessible
quotient, and achieved in the case of Schneider's equivalence
in each localization of the cover. In a sense, this
is a globalization of Schneider's equivalence, i.e.
a special case of the
(effective) descent method for noncommutative nonaffine base (quotient) space.

One could naively always {\it define}
$Qcoh(E/B)$ as being simply ${}_{\cal E}{\cal M}^{\cal B}$. However,
this may not be always achieved, if we want to equip this category with 
local coordinate charts. 
The ${\cal B}$-comodule algebra ${\cal E}$ itself (put coaction aside
for a moment) is not only defining the category of ${\cal E}$-modules, 
it is itself also a choice of the underlying ring of that category
within its Morita equivalence class.
A global version of that choice is achieved as follows. 
\newline ${\cal E}$ is
viewed as a {\it relative} noncommutative scheme~\cite{Ros:NcSch}
over the category of modules over $\genfd$. 
The inverse image of the distinguished object 
in the category of $\genfd$-modules is the structure sheaf ${\cal O}_{\cal E}$.
For any perfect localization $Q$, the inverse image $Q^* ({\cal O}_{\cal E})$
is a $\genfd$-algebra. One may require that the noncommutative quotient scheme 
$E/B$ over $\genfd$ also comes equipped with ${\cal O}_{E/B}$. 
Realizing this structure
sheaf for some category close to ${}_{\cal E}{\cal M}^{\cal B}$ is the
basic meaning of the {\bf construction of the quotient space}. 
Our approach is to localize enough to get 
${\cal O}_{E/B}$ from local coinvariants,
or at least to approximate it. ``Size'' of an ${\frak U}$-accessible quotient
is a naive measure of the approximation.
An ideal case is when each localization of ${\cal E}$ 
is a Hopf-Galois extension:
intuitively, the local coinvariants then patch 
to the full structure sheaf. 
One can extend the theory beyond the affine case.  
In order to make sense of coactions, 
(nonaffine) noncommutative $\genfd$-scheme 
${\rm Spec}({\cal E})$ may be required to come with an affine cover
by compatible perfect localizations.  Those are given coactions (on the
pieces of ${\cal O}_{\cal E}$). The coactions should
agree when viewed at the level of categories of modules
(in the commutative case only some $B$-schemes are of this form, but it is
strictly larger than affine $B$-schemes). 
To achieve further progress 
it is promising to rely on simplicial and cohomological methods;
the ideas of cohomological descent~\cite{Del:HodgeIII}
and the noncommutative 
\v{C}ech cohomologies~\cite{Ros:NcSch, vOystWill:cohschematic}
provide a framework. Cf. also~\cite{Beck,BernLunts,
Del:tan,Giraud,Jard:simp,LuntsRos2,Moerd:class,Moerd:clpap,Moerd:clpap}.
Very recently, {\sc V.~Lunts} suggested that it might be 
useful to consider more general flat covers
and resolutions of would-be quotient space
than the covers coming from localizations,
and to use the corresponding (co)simplicial objects
to define and investigate quotients of comodule algebras 
in full generality~\cite{LuntsSkoda}
(even nontrivial 
$\rho$-compatible localization covers are not always available). 
However, in the rest of the paper, I will
present my earlier work focusing on localization 
approach, which for some purposes suffices.

\section*{10. Quantum bundles.}
A right (left) Zariski 
{\bf locally trivial principal ${\cal B}$-bundle}
is a right (left) ${\cal B}$-comodule algebra $({\cal E},\rho)$
together with a right (left) Zariski local trivialization. 
A {\bf local} right (left) {\bf trivialization} of 
$({\cal E},\rho)$ consists of

$\,\,\,\,\,\,\,\bullet\,\,$ 
a {\bf cover} of ${\cal E}$
by {\bf $\rho$-compatible perfect localizations} 
$i_\lambda : {\cal E} \rightarrow {\cal E}_\lambda$,

$\,\,\,\,\,\,\,\bullet\,\,$
a family $\{ \gamma_\lambda : {\cal B} \rightarrow {\cal E}_\lambda\}$
of right (left) ${\cal B}$-comodule algebra homomorphisms.

Here the ${\cal B}$-comodule structure on ${\cal E}_\lambda$ 
is the one induced by $\rho$-compatibility.
\vskip .03in
\def\ooe{{\cal E}}

Local triviality implies that each 
$({\cal E}_\lambda,i_\lambda,\gamma_\lambda)$
is in the category of higher smash products ${\cal C}(U,{\cal B})$
for some ${\cal U}$ (trivial quantum principal bundle).

Let now $(M,\rho_M)$ be a {\it left} ${\cal B}$-comodule. Define
\begin{equation}\label{eq:kappaibar}\begin{array}{l}
\kappa_\lambda = (\mu \otimes {\rm id})
({\rm id} \otimes \gamma_\lambda \otimes 
{\rm id})({\rm id} \otimes \rho_M)
: {\cal E}_\lambda^{{\rm co}\cal B} \otimes M
\rightarrow {\cal E}_\lambda \otimes M,\\
\bar{\kappa}_\lambda = (\mu \otimes {\rm id})
({\rm id} \otimes (\gamma_\lambda\circ S) \otimes {\rm id})
({\rm id} \otimes \rho_M) : {\cal E}_\lambda \otimes M \rightarrow 
{\cal E}_\lambda\otimes M.
\end{array}\end{equation}
Equivalently, $\kappa_\lambda$ is
the unique map such that 
\[ \kappa_\lambda (f\otimes m) = 
\sum f\,\gamma_\lambda(m_{(-1)}) \otimes m_{(0)}.\]
These data are enough to define the {\bf associated (vector) bundle} $\xi_M$,
when the cover $\{\lambda \in \Lambda\}$ is finite.
Denote
$\Gamma_\lambda\xi_M = {\cal E}_\lambda^{{\rm co}\cal B} \otimes M$
and define $\Gamma_\lambda\xi_M$ to be a $\genfd$-space 
of local sections of associated bundle over localization $\lambda$.
For other localizations $Q_\mu$ (not included in the trivializing
cover), a version of the globalization lemma 
(Theorem 6.2 in~\cite{Rosen:book}) and elementary arguments justify 
the definition by a ``$\kappa$-twisted descent'' formula
\[ \Gamma_\mu \xi_M 
= \{ 
F = \prod_{\lambda \in \Lambda} F_\lambda \in 
\prod_{\lambda \in \Lambda} Q_\mu(\Gamma_\lambda \xi_M) 
,\,
\left\lbrace\begin{array}{l}
(i^{\lambda\mu}_{\lambda\mu,\lambda'\mu'} \otimes \id)
\kappa_\lambda F_\lambda
= (i^{\lambda'\mu'}_{\lambda\mu,\lambda'\mu'} \otimes \id)\kappa_{\lambda'} 
F_{\lambda'}\\
(i^{\lambda\mu}_{\lambda'\mu',\lambda\mu} \otimes \id)\kappa_\lambda F_\lambda
= (i^{\lambda'\mu'}_{\lambda'\mu',\lambda\mu} \otimes \id)\kappa_{\lambda'} 
F_{\lambda'}
\end{array}\right.
\}
\]
This mentioned, below we limit ourselves to the original
cover and the global sections.

\th{Proposition}{1.} { 
$\bar{\kappa}_\lambda \circ \kappa_\lambda = \id$.
In particular,
$\kappa_\lambda$ is an isomorphism onto its image and 
$\kappa_\lambda \circ \bar{\kappa}_\lambda 
|_{{\rm Im} \kappa_\lambda} = \id$. 
}

The {\bf cotensor product} of 
a left ${\cal B}$-comodule M and 
right ${\cal B}$-comodule ${\cal E}$ is defined by
\[ {\cal E} \Box_{{\cal B}} M = \{ H \in {\cal E}\otimes M \,|\, 
({\rm id} \otimes \rho_M) H  = (\rho_E \otimes {\rm id}) H \}. \]

The cotensor product has {\it a priori} only a structure of
a $\genfd$-vector space. If ${\cal E}$ is a ${\cal B}$-comodule algebra
which is the total space of a Hopf-algebraic trivial principal bundle
(which we identify with the smash product
${\cal E}^{{\rm co}{{\cal B}}} \sharp {\cal B}$), then 
it is elementary that ${\cal E}^{{\rm co}{{\cal B}}} \otimes M$
is canonically isomorphic to the cotensor product as a $\genfd$-vector space.
We extend this result to the setting of 
locally trivial bundles.

\th{Lemma}{2.~\cite{Skoda:thes,Skoda:loc-coinv1}}{ 
Image of ${\cal E} \Box_{\cal B} M$ in localization
${\cal E}_\lambda \otimes M$ is contained in
${\rm Im}\,\kappa_\lambda$.}

The algebra of global sections of {\bf associated bundle} $\Gamma \xi_M$
is (in this setup) given by
\[ \Gamma \xi_M = \{ F = \prod_{\lambda \in \Lambda} F_\lambda
\in \prod_{\lambda \in \Lambda} \Gamma_\lambda \xi_M,\,
\left\lbrace\begin{array}{l}
(i^\lambda_{\lambda,\lambda'} \otimes \id)\kappa_\lambda F_\lambda
= (i^{\lambda'}_{\lambda,\lambda'} \otimes \id)\kappa_{\lambda'} 
F_{\lambda'}\\
(i^\lambda_{\lambda',\lambda} \otimes \id)\kappa_\lambda F_\lambda
= (i^{\lambda'}_{\lambda',\lambda} \otimes \id)\kappa_{\lambda'} 
F_{\lambda'}
\end{array}\right.\}. \]
\th{Theorem} {3.}{ 
$\Gamma \xi_M$ is naturally isomorphic to ${\cal E}\Box_{\cal B} M$
as a $\genfd$-vector space.
}

Proof. By the globalization lemma, there is an isomorphism 
$i_\Lambda$ from ${\cal E}$ to the equalizer of the diagram
$\prod {\cal E}_\lambda \rightrightarrows \prod {\cal E}_{\lambda,\lambda'}$.
Tensoring $i_\Lambda$ and the 
equalizer diagram with ${\rm id}_M$ over the ground field
$\genfd$ yields an isomorphism $i_\Lambda \otimes {\rm id}_M$ from
${\cal E}\otimes M$ to the equalizer of
$\prod {\cal E}_\lambda \otimes M
\rightrightarrows \prod {\cal E}_{\lambda,\lambda'}\otimes M$.
A direct check shows that $\prod\rho_{{\cal E}_\lambda}\otimes {\rm id}_M$
and $\prod {\rm id}_{{\cal E}_\lambda}\otimes \rho_M$ agree on
the image of map $\bar{\kappa}_\lambda
(i_\lambda \otimes \id_M)$ defined on 
$\prod_\lambda {\cal E}_\lambda^{{\rm co}{\cal B}}\otimes M$. 
As our localization endofunctors are flat, 
and the coactions $\rho_{{\cal E}_\lambda}$ extend $\rho_{\cal E}$,
the covering property implies that $\prod \rho_{\cal E}\otimes {\rm id}_M$
agrees with $\prod {\rm id}_{{\cal E}_\lambda}\otimes \rho_M$
in the equalizer.
The isomorphism $\Gamma \xi_M \rightarrow {\cal E}\Box_{\cal B} M$
is well defined by 
$(i_\Lambda \otimes \id_M)^{-1} \prod_\lambda \kappa_\lambda$.
Namely, its inverse is $H \mapsto \prod_\lambda \bar{\kappa}_\lambda
(i_\lambda \otimes \id_M)H$.
This follows by Lemma 2 and the properties of $\kappa_\lambda$-s
listed above. $\Box$
 
\th{Corollary}{4.} {Let ${\cal G}$ be an Hopf algebra,
${\cal B}$ a quantum subgroup (quotient Hopf algebra). 
If $({\cal G}_\lambda,\rho_{\cal B})$
is a Zariski locally trivial principal $\rho_{\cal B}$-bundle, 
then the inducing functor for comodules
${\rm Ind}^{\cal G}_{\cal B}\,: {\cal M}^{\cal B}
\rightarrow {\cal M}^{\cal G}$ can be realized by gluing local sections.}

Here $\pi : {\cal G}\rightarrow {\cal B}$ is onto map of Hopf
algebras. The induced comodule can be described as
${\rm Ind}^{\cal G}_{\cal B} M = {\cal G}\Box_{\cal B} M$
where $\rho_{\cal B} : {\cal G} \rightarrow {\cal G}\otimes {\cal B}$
is given by $\rho_{\cal B} = (\id \otimes \pi)\Delta_{\cal G}$.
\vskip .022in

One can interpret the global sections functor $\Gamma$
as an appropriate $0$--th~\v{C}ech--type cohomology~$\check{H}^0$. 
Similarily, the application of the methods in~\cite{Ros:NcSch}
imply a generalization of Corollary~4: the composition
$\check{H}^i \circ \xi_{\_}$, where $\check{H}^i$ are the appropriate higher
\v{C}ech-type cohomologies and $\xi_{\_}$ is the associated bundle
functor $M \mapsto \xi_M$,
realize the higher derived functors $R^i {\rm Ind}^G_B$ 
of the induction functor for comodules.
\vskip .018in

We use (co)matrix notation for coaction
to introduce the transition matrices for $\xi_M$. 
Supose $\{m_\beta\}$ is a basis of $M$ and the ${\cal B}$-coaction is given by 
\[ \rho \,:\, m_\beta \mapsto \sum_\alpha m^\alpha_\beta \otimes m_\alpha \]
Represent element $F_\lambda \in \Gamma_\lambda \xi_M$ as
\[ F_\lambda = \sum_\alpha f^\alpha_\lambda \otimes m_\alpha.\]
The condition $\kappa_\lambda F_\lambda = \kappa_{\lambda'} F_{\lambda'}$
(in both consecutive localizations) reads
\[ \sum_{\alpha, \beta} f^\alpha_\lambda \gamma_\lambda(m^\beta_\alpha)
\otimes m_\beta = 
\sum_{\alpha, \beta} f^\alpha_{\lambda'} \gamma_{\lambda'}(m^\beta_\alpha)
\otimes m_\beta. \]
As $\{m_\beta\}$ is a basis, this is equivalent 
to the statement
\[ (\forall \beta)\,\,
\sum_{\alpha} f^\alpha_\lambda \gamma_\lambda(m^\beta_\alpha) = 
\sum_{\alpha} f^\alpha_{\lambda'} \gamma_{\lambda'}(m^\beta_\alpha). \]
Multiplying this identity by $\gamma_\lambda(Sm^\gamma_{\beta})$
and summing over $\beta$ we obtain
\begin{equation}\label{eq:fortranmatrix}
 \sum_{\alpha, \beta} f^\alpha_\lambda \gamma_\lambda(m^\beta_\alpha)
\gamma_\lambda(Sm^\gamma_\beta) = 
\sum_{\alpha, \beta} f^\alpha_{\lambda'} \gamma_{\lambda'}(m^\beta_\alpha)
\gamma_\lambda(Sm^\gamma_\beta).
\end{equation}
As $\gamma_\lambda$ is an algebra map,
$\sum_\alpha \gamma_\lambda(m^\beta_\alpha)
\gamma_\lambda(Sm^\gamma_\beta) = \delta^\gamma_\alpha 1$
what simplifies the left--hand side to $f^\gamma_\lambda$. Denote
\begin{equation}\label{eq:tranmatrix}
 \left(\frak{M}_{\lambda',\lambda}\right)^\gamma_\alpha 
\, := \sum_\beta \gamma_{\lambda'}(m^\beta_\alpha)
\gamma_\lambda(Sm^\gamma_\beta) \,\,\,\,\in \, {\cal E}_{\lambda\vee\lambda'}.
\end{equation}
This defines a matrix ${\frak M}_{\lambda,\lambda'}$ 
which is called the {\bf transition matrix between charts $\lambda$ and
$\lambda'$}. 
In this notation~(\ref{eq:fortranmatrix}) becomes
\[ f^\gamma_\lambda = \sum_\alpha f^\alpha_{\lambda'} 
\left({\frak M}_{\lambda',\lambda}\right)^\gamma_\alpha. \]
It follows from definition~(\ref{eq:tranmatrix}) that
\[ {\frak M}_{\lambda,\lambda} = I,\,\,\,\,\
{\frak M}_{\lambda,\mu} {\frak M}_{\mu,\nu} = 
{\frak M}_{\lambda,\nu}, \]
where the upper index of matrix $\frak{M}_{\lambda',\lambda}$ 
is the column index, rather than the row index 
(the latter would be along conventions for other matrices in this work).

Finally, we comment on the restriction of 
(algebras of global sections of) bundles
to Ore localizations. 
\th{Theorem}{5.~\cite{Skoda:thes, Skoda:loc-coinv1}}{
 There exist a Hopf algebra ${\cal B}$, 
and a ${\cal B}$-comodule algebra $({\cal E},\rho)$ such that
\begin{itemize}
\item The algebra of (global) coinvariants ${\cal E}^{{\rm co}{\cal B}}$
is commutative (possible to choose the polynomial ring $\genfd[u]$).
\item ${\cal E}$ is isomorphic to a smash product of ${\cal B}$
and ${\cal E}^{{\rm co}{\cal B}}$
with coaction $\it \otimes \Delta$.
\item There exist a multiplicative set 
$U \subset {\cal E}^{{\rm co}\cal B}$
(automatically Ore by commutativity),
such that there is no pair of a $\rho$-compatible left Ore set
$S \subset {\cal E}$ 
and an isomorphism $\zeta$ of ${\cal B}$-comodule algebras
\[ \zeta : (S^{-1}{\cal E})^{{\rm co}{\cal B}} \rightarrow 
U^{-1}({\cal E}^{{\rm co}\cal B}),\]
such that $\zeta \circ i_S |_{{\cal E}^{{\rm co}\cal B}} = i_U$.
\end{itemize}
}
This phenomenon is nonexistent in commutative and
generic in noncommutative case, even if, as here, all the global coinvariants
pairwise commute. In~\cite{SchMon:quotients} 
some conditions preventing such phenomena
were considered for certain kinds of localizations using filters of
ideals. In~\cite{LuntsRos2} graded actions and actions by differential
operators were considered and, in their setup, guaranteed that
module algebra actions extend to localizations. In the language
of modules and dimodules, the simplest reasoning about these
phenomena is as follows. Coaction on the whole ring does not
necessarily extend to the localized ring as we know. If they do,
then the localization induces functor $Q^{\cal B}$ on 
${}_{\cal E}{\cal M}^{\cal B}$, so if the latter category is equivalent
to ${}_{\cal U}{\cal M}$ (say, in faithfully flat Hopf-Galois case),
we obtain an induced localization in the base for free. 
The induced radical $\sigma_U$ in 
the base category ${}_{\cal U}{\cal M}$ is given by $M \mapsto N$ whenever 
$\sigma^{\cal B} : M \otimes {\cal B} \mapsto N \otimes {\cal B}$.
One can not apriori say if the induced localization of ${}_{\cal U}{\cal M}$
is of some specific type. Say, if we start with a
$\rho$-compatible Ore localization
of ${}_{\cal E}{\cal M}$,
is the induced localization of  ${}_{\cal U}{\cal M}$
Ore or only a perfect localization ? This is not very tractable
question. 
A ${\cal B}$-comodule algebra $({\cal E},\rho)$
is {\bf sufficiently localized} with respect to
the family of all $\rho$-compatible left Ore localizations 
if for every such $S$ there is a left Ore set $U$ 
in ${\cal U} = {\cal E}^{{\rm co}{\cal B}}$ such that 
\[  i_S|{\cal E}^{{\cal B}} :
{\cal E}^{{\cal B}} \rightarrow (S^{-1}{\cal E})^{\cal B}\]
is canonically isomorphic to an Ore localization of ${\cal E}^{\cal B}$
with respect to some Ore set $U \subset {\cal E}^{\cal B}$.
It is interesting to compare this notion with other mentioned conditions
on the relation between the base and total algebra. 
For example, by the fundamental theorem of Hopf modules, 
each $\rho$-compatible flat localization of ${}_{\cal B}{\cal M}^{\cal B}$
is equivalent to a flat localization of the category of $\genfd$-vector
spaces. As all such localizations are trivial, 
{\it it is not possible to find a nontrivial Ore
localization ${\cal G}[T^{-1}]$ 
of a Hopf algebra ${\cal G}$
such that the comultiplication $\Delta_{\cal G}$ extends multiplicatively to
a coaction $\Delta_T$ in the localization}. 
This is an extension of the basic fact that
any $G$-invariant subset of an algebraic group $G$ is $G$ itself. 

\vskip .022in
Start with a localization $Q_{\cal U}$ of the base category
${}_{\cal U}{\cal M}$,
and suppose it is {\it abstractly} equivalent to the category
of dimodules ${}_{\cal E}{\cal M}^{\cal B}$ 
on the total algebra ${\cal E}$. Then, of course, we have 
a localization $Q^{\cal B}$ on that category, 
but there is no {\it apriori} reason to
believe that $Q^{\cal B}$ is induced by any localization $Q$ defined 
on the ambient category ${\cal M}$. For example, the same
module $M$ can have two different Hopf module structures $M_1$ and
$M_2$, and $Q^{\cal B}(M_1)$ may disagree with $Q^{\cal B}(M_2)$
as an ${\cal E}$-module. 
If the abstract category equivalence is however the
natural one -- as in the Schneider's
theorem, that is via $()^{{\rm co}{\cal B}}$ and $\otimes_\genfd {\cal B}$
-- then this bad phenomenon is prevented.

Ore sets in the coinvariants algebra ${\cal U}$ (for trivial bundle) 
may correspond to more general (than Ore) 
$\rho$-compatible localizations in the total space, so taking
just the ``topology'' generated by Ore localizations would
mean that the geometric projection is not continuous in the obvious sense,
even in otherwise well--behaved noncommutative examples. 

\section*{11. Cohn localization.}
~\cite{Cohn:free, Cohn:invloc} 
Let $R$ be a possibly noncommutative ring, and
$\Sigma$ a given set of square {\it matrices} of
possibly different sizes with entries in $R$.
Map $f: R \rightarrow S$ of rings is $\Sigma$-inverting
if each matrix in $\Sigma$ is mapped to an invertible
matrix over $S$.  A $\Sigma$-inverting ring map 
$i_\Sigma : R \rightarrow R_\Sigma$ 
is called {\bf Cohn localization} (or universal
$\Sigma$-inverting localization) if for
every $\Sigma$-inverting ring map $f: R \rightarrow S$
there exist a unique ring map $\tilde{f} : R_\Sigma \rightarrow S$
such that $f = \tilde{f}\circ i_\Sigma$.

A set $\Sigma$ of matrices is called {\bf multiplicative}
if $1 \in \Sigma$ and, for any $A,B \in \Sigma$ and matrix $C$
of right size over $R$, the matrix 
$\left(\begin{array}{cc}A & C \\ 0 & B \end{array}\right) \in \Sigma$.
If $\Sigma$ is multiplicative, then the elements of $R_\Sigma$ can be
obtained as the components of the solutions to all the matrix 
equations $Au = a$, where $a$ is a column over $R$ and $A \in \Sigma$.
To any multiplicative $\Sigma$, {\sc P.~M. Cohn}~\cite{Cohn:invloc}
associates an idempotent preradical $\sigma_\Sigma$.
Its left-module variant is given by
\[ M \mapsto \sigma_\Sigma(M) 
= \{ m \in M \,|\, \exists u = (u_1,\ldots, u_n)^T, 
\, \exists i, \, m = u_i\mbox{ and } \exists A \in \Sigma,
\,\, A u = 0\}. \]
Associated torsion theory is not always hereditary 
(i.e. a submodule of $\sigma_\Sigma$-torsion
module is not necessarily $\sigma_\Sigma$-torsion).
Equivalently, $\sigma_\Sigma$ is not necessarily left exact.
Hereditary torsion theories correspond to Gabriel localizations.
In practice, the main problem with Cohn localization
is that it is usually hard to determine the 
kernel of the localization map.

\section*{12. Quasideterminants and noncommutative Gauss decomposition.}

Inverting matrices with noncommutative entries
plays role not only in Cohn localization but also generally 
in solving linear systems of equations in noncommutative variables. 
The main tool to do calculations with
such matrix inverses and perform Gauss-type decompositions are
{\it quasideterminants} of Gel'fand and Retakh
~\cite{GelRet:ncr,KrobLec,GelRet:graph,GelRet:quasiI,Skoda:thes}.

Let $A\in M_n(R)$ be a $n\times n$ matrix
over an arbitrary noncommutative unital associative ring $R$. 
Suppose rows and columns of $A$ are labeled. Let us
choose a row label $i$ and a column label $j$.
By $A^{\hat{i}}_{\hat{j}}$ we'll denote the $(n-1)\times(n-1)$ matrix obtained
from $A$ by removing the $i$-th row and the $j$-th column.
The $(i,j)$-th {\bf quasideterminant} $|A|_{ij}$ is
\[\fbox{$ |A|_{ij} = a^i_j - 
\sum_{k \neq i, l\neq j} a^i_l  (A^{\hat{i}}_{\hat{j}})^{-1}_{lk} a^k_j
$}\]
provided the right-hand side is defined.

At most $n^2$ quasideterminants of a given $A \in M_n(R)$
may be defined. If all the $n^2$ quasideterminants $|A|_{ij}$ exist 
and are invertible then the inverse $A^{-1}$ of $A$ exist in 
$A \in M_n(R)$ and
\begin{equation}\label{eq:quasiInv}\index{$|A|_{ji}$} 
(|A|_{ji})^{-1} = (A^{-1})^i_j. 
\end{equation}

Suppose now we are given an equation of the form
\[ A x = \xi \]
Define thus $A(j,\xi)$ as the $n \times n$ matrix whose entries are
the same as of $A$ except that the $j$-th column is replaced by $\xi$.
Then the noncommutative {\bf left Cramer's rule} says
\[\fbox{$ |A|_{ij} x^j = |A(j,\xi)|_{ij} $} \]

\vskip .012in
Quasideterminants for classical matrices are up to sign
ratios of a determinant of a matrix and the determinant
of $(n-1)\times (n-1)$ submatrix, as one can see by 
remembering the formula for matrix inverse in terms of
cofactor matrices.

\vskip .012in
Quasideterminants satisfy a number of useful identities,
namely, analogues of classical Muir's law of extensionality, 
Jacobi inversion formula, Laplace expansion formulas etc. 
There is also a new property called heredity.
Descriptively it is a compatibility of taking quasideterminant
with partitioning of a matrix (into square matrix blocks).
In other words, one can take quasideterminants of
partitioned matrices in stages.

\vskip .012in
{\bf Noncommutative Gauss decomposition} 
of a matrix $G$ is the decomposition
\begin{equation}\label{eq:ncGauss} G = UA \end{equation}
here $A$ is a lower triangular matrix
(with possibly noncommutative entries) and $U$ is an upper unidiagonal matrix
(i.e. $u^i_j = 0$ for $i>j$, and $u^i_j = 1$ for $i=j$).

\vskip .002in
The problem~(\ref{eq:ncGauss}) is equivalent to the set of $n^2$ equations
\[ g^i_j = \left\lbrace \begin{array}{ll} a^i_j + 
        \sum_{k > i \geq j} u^i_k a^k_j
         \,,\,\,\,\,\, i \geq j, \\
         \sum_{k \geq j > i} u^i_k a^k_j
         \,,\,\,\,\,\, i < j,
\end{array} \right.\]
In terms of quasideterminants the solution is
\begin{equation}\label{eq:sol-Gauss}
\begin{array}{l}
  a^i_j = \left|  \begin{array}{cccc}
        g^i_j & g^i_{i+1} & \cdots & g^i_n \\
g^{i+1}_j & g^{i+1}_{i+1} & \cdots & g^{i+1}_n \\
        \ldots & \ldots & \ldots & \ldots\\
        g^n_j & g^n_{i+1} & \cdots & g^n_n
\end{array}     
        \right|_{ij},\,\,\,\,\,\,\, i \geq j,
\\
  u^i_j = \left| \begin{array}{cccc}
        g^i_j & g^i_{j+1} & \cdots & g^i_n \\
g^{j+1}_j & g^{j+1}_{j+1} & \cdots & g^{j+1}_n \\
        \ldots & \ldots & \ldots & \ldots\\
        g^n_j & g^n_{j+1} & \cdots & g^n_n
\end{array} \right|_{ij}
        \left| \begin{array}{cccc}
        g^j_j & g^j_{j+1} & \cdots & g^j_n \\
g^{j+1}_j & g^{j+1}_{j+1} & \cdots & g^{j+1}_n \\
        \ldots & \ldots & \ldots & \ldots\\
        g^n_j & g^n_{j+1} & \cdots & g^n_n
\end{array}
        \right|^{-1}_{jj},\,\,\,\,\,\, i < j
\end{array}
\end{equation}
whenever all the principal (=lower right corner) 
quasiminors (=quasideterminants of submatrices) 
exist and they are invertible. We suggest to the reader 
to picture (positions of) rows and columns in $G$ which are
involved in the submatrices of $G$ in~(\ref{eq:sol-Gauss}).


\section*{13. Matrix bialgebras and Hopf algebras.} 
{\bf 13.1} The full ring  $M(n,\genfd)$ 
of $n\times n$ matrices with (commutative) entries in a field $\genfd$ is
isomorphic to $\genfd^{n^2}$ as a $\genfd$-vector space. 
This isomorphism induces a structure of affine $\genfd$-variety
on $M(n,\genfd)$.
The regular functions on that variety are polynomials in matrix entries.
Introduce $n^2$ regular functions
\[ t^i_j : M(n,\genfd) \rightarrow \genfd,\,\,\,\,\,\,
\,\,\,\,\,t^i_j(a) = a^i_j, \,\,\,\,\,\,
a \in M(n,\genfd),\,\,\, i,j= 1,\ldots n. \]

Then ${\cal O}(M(n,\genfd)) \cong \genfd [t^1_1, t^1_2, \ldots, t^n_n]$ 
is the ring of global regular functions on $M(n,\genfd)$. 

\vskip .03in
Let ${\cal G}$ be a bialgebra, possibly noncommutative,
over a field $\genfd$ and
$G = (g^i_j)^{i = 1,\ldots, n}_{j = 1,\ldots, n}$ 
an $n\times n$-matrix over ${\cal G}$.
${\cal G}$ is a {\bf matrix bialgebra} with basis $G$ if
the set of entries of $G$ generates ${\cal G}$ and
if the comultiplication $\Delta$ and counit $\epsilon$ satisfy
\[\begin{array}{lll} \Delta G = G \otimes G\,&\mbox{ i.e. }\,& 
\Delta g^i_j = \sum_{k = 1}^n g^i_k \otimes g^k_j\\
\epsilon G = {\bf 1}\,&\mbox{ i.e. }\,&\epsilon (g^i_j) = \delta^i_j
\end{array}\]
The free (noncommutative) associative algebra $F$ 
on $n^2$ generators $f^i_j$
has a unique coalgebra structure making it a
matrix bialgebra with basis $(f^i_j)^{i = 1,\ldots, n}_{j = 1,\ldots, n}$. 
We call it the {\bf 
{free matrix bialgebra of rank $n^2$}}.
Every bialgebra quotient of that bialgebra is a matrix bialgebra.

\vskip .03in
{\bf 13.2} {\it Free Hopf algebras.}
We are going to sketch the construction of free Hopf algebras
due {\sc Takeuchi}~\cite{freeHopf}.
Let $C$ be a coalgebra over $\genfd$. Let $C_i = C$
for $i$ even nonnegative integer, and $C_i = C_i^{{\rm cop}}$
(the cooposite coalgebra of $C$) for odd positive integer. Then define
$V$ to be the external direct sum (coproduct) of coalgebras $C_i$,
\[ V = \coprod C_i. \]
The tensor algebra $T(V)$ of $V$, as of any other coalgebra,
has a unique bialgebra structure such that the natural inclusion
$i_V : V \rightarrow T(V)$ is a coalgebra map. It holds
$T(V^{\rm cop}) \cong T(V)^{\rm cop}$.
Then define $\genfd$--linear map
\[ S_V : V \rightarrow V^{{\rm cop}}\,\,\,
\mbox{ by } \,\,(v_0, v_1, v_2,\ldots) = 
(0, v_0, v_1, v_2, \ldots). \]
There is a unique bialgebra map 
$S : T(V) \rightarrow T(V)^{\rm cop}$ extending $S_V$.
Let $I_S$ be the 2-sided ideal in $T(V)$
generated by all elements of the form
\[ \sum c_{(1)} S(c_{(2)}) - \epsilon(c) 1\,\,\mbox{ and }
\sum S(c_{(1)}) c_{(2)} - \epsilon(c) 1,\,\,c \in C_i, \,i = 1,2, \ldots \]
This 2-sided ideal is a biideal and $S(I_S) \subset I_S$, hence
it induces a bialgebra map
\[ S : T(V)/I_S \rightarrow (T(V)/I_S)^{\rm cop}. \]
It follows that $H(C) = T(V)/I_S$ is a Hopf algebra, the {\bf free
Hopf algebra} on $C$.
For any Hopf algebra $H'$ and a coalgebra map
$\phi : C \rightarrow H'$ there is a unique Hopf algebra
map $\phi' : H(C) \rightarrow H'$ such that
$\phi' \circ i = \phi$ where $i : C \rightarrow H(C)$ is the 
composition of inclusion 
into $T(V)$ and projection $T(V) \rightarrow H(C)$.

\vskip .013in
Notice that, for any family $\{ C_i \}$ of coalgebras,
$T(\coprod_i C_i) = \coprod_i T(C_i)$ where
the coproduct on the RHS is the coproduct in the category of
bialgebras, what is the free product of algebras with natural
induced coalgebra structure. In the case when the index set are
nonnegative integers and $C_{i+1} = C_i^{\rm cop}$, 
the LHS specializes to an intermidiate
stage in building $H(C)$, so it is a generalization of it. 

\vskip .011in
{\sc Manin}~\cite{ManinQGNG} generalized the RHS. 
He replaces $T(C_i)$ by any  bialgebra $B_i$ 
with $B_{i+1} = B_i^{\rm cop, op}$.
Notice that the algebra structure is also opposite between
even and odd cases (superfluous/unvisible condition in the case of $T(C_i)$).
Let ${\bf B} = \coprod_i B_i$ and $S : {\bf B} \rightarrow {\bf B}$
be again defined by a shift in index by $+1$. Then the 2-sided
ideal $I_S \subset {\bf B}$ generated by 
\[ \sum b_{(1)} S(b_{(2)}) - \epsilon(b) 1\,\,\mbox{ and }
\sum S(b_{(1)}) b_{(2)} - \epsilon(b) 1,\,\,b \in B_i, \]
is $S$-stable and 
$H(B) = {\bf B}/I_S$ is a Hopf algebra,
the {\bf Hopf envelope} of $B$.

It satisfies the following universal property:
for any Hopf algebra $H'$ and a bialgebra map
$\phi : B_0 \rightarrow H'$ there is a unique Hopf algebra
map $H(\phi) : H(B) \rightarrow H'$ such that
$H(\phi) \circ i = \phi$ where $i : B_0 \rightarrow H(B)$ is the 
composition of inclusion into ${\bf B}$ and projection 
${\bf B} \rightarrow H(B)$.

\vskip .026in
{\bf 13.3} {\it Matrix Hopf algebras.}
One can specialize the construction of Hopf envelope
to any matrix bialgebra on basis 
$T = (t^i_j)^{i = 1,\ldots,n}_{j = 1,\ldots,n}$.
A {\bf matrix Hopf algebra} ${\cal G}$ with basis $T = (t^i_j)$
is a Hopf algebra which possess a matrix subbialgebra
$B$ with basis $T$ such that the map $H(\id):H(B)\rightarrow{\cal G}$
is onto. $B$ is not necessarily a matrix bialgebra with respect to
that basis: e.g. commutative
coordinate ring of $GL(n,\genfd)$ is not a matrix bialgebra 
with respect to the obvious basis $T$, as one needs
to introduce inverse of the determinant. That can be, of course, repaired
by enlarging the basis. On the other hand, ${\cal O}(SL(n,\genfd))$ is
a matrix bialgebra and a matrix Hopf algebra with 
the same standard basis $T$.

\vskip .026in
{\bf Free matrix Hopf algebra} ${\cal NGL}(n,\genfd)$
is the free Hopf algebra generated by the free matrix coalgebra
${\cal M}_C(n,\genfd)$ (the latter is just the
$\genfd$-vector space $C$ of dimension $n^2$ with matrix comultiplication).
Equivalently, ${\cal NGL}(n,\genfd)$ is the Hopf envelope
of the free matrix bialgebra ${\cal NM}(n,\genfd)$ on $n^2$ generators.
Notice that ${\cal M}_C(n,\genfd)$ is finite dimensional, 
${\cal NM}(n,\genfd)$ is finitely generated and
${\cal NGL}(n,\genfd)$ is apriori not even finitely generated.

\th{Proposition}{6.}
{
Let $\Sigma_1 \subset \Sigma_2$ be two
sets of submatrices of 
the basis $T$ of ${\cal NM}_n = {\cal NM}(n,\genfd)$.
Then the natural maps
$\Sigma_1^{-1} {\cal NM}_n \rightarrow \Sigma_2^{-1} {\cal NM}_n$
are 1-1.
}

Proof. It is known~\cite{Cohn:free} that there is an algebra embedding of
any {\it free} associative algebra into a skewfield ${\cal K}$.
Let $j : {\cal NM}_n\rightarrow {\cal K}$ be such an embedding.
By the universal property of the localization there is a unique map
$j_\Sigma : \Sigma^{-1} {\cal NM}_n \rightarrow {\cal K}$ provided
$\Sigma$ consists of matrices which are invertible over ${\cal K}$.
For that it is sufficient to prove that all the quasideterminants
exist and they are invertible (nonzero). 
Quasideterminants can be defined inductively by size,
provided all the involved quasiminors at every step are nonzero.
In our case we deal with quasiminors of $T$ only. 
Their images under the projection 
${\cal K} \rightarrow {\cal K}/[{\cal K},{\cal K}]$
onto the commutative field ${\cal K}/[{\cal K},{\cal K}]$ are
apparently nonzero, hence they are nonzero over ${\cal K}$, 
hence invertible.

Using the universal property again,
we see that the chain of natural maps 
${\cal NM}_n \rightarrow \Sigma_1^{-1}{\cal NM}_n \rightarrow 
\Sigma_2^{-1}{\cal NM}_n\rightarrow {\cal K}$ composes to $j$, hence
all the maps in the chain are 1-1. $\Box$

\vskip .03in
{\bf Free triangular matrix Hopf algebra} ${\cal NB}(n,\genfd)$
can be obtained directly in analogy
with ${\cal NGL}(n,\genfd)$. The only difference is that
the coalgebra (bialgebra) generators $t^i_j$ for 
$i<j$ of $C = {\cal M}_c(n,\genfd)$ (${\cal NM}(n,\genfd)$
respectively) are set to zero before taking the construction. 
Alternatively,
one can start with ${\cal NGL}(n,\genfd)$ and
consider the 2-sided ideal generated by $(t^i_j)_r$ for 
all $r = 0,1,2,\ldots$ and all $i<j$. This biideal is 
obviously $S$-stable, hence induces a quotient map of Hopf algebras
$\pi : {\cal NGL}(n,\genfd) \rightarrow {\cal NB}(n,\genfd)$
whose image is isomorphic to free triangular matrix Hopf algebra.
In the latter setup, we call ${\cal NB}(n,\genfd)$, or, more
precisely, pair $({\cal NB}(n,\genfd),\pi)$, the
{\bf noncommutative Borel subgroup} in ${\cal NGL}(n,\genfd)$.

\vskip .01in
Given a permutation $\sigma$ on $n$ letters, 
and a $n \times n$ matrix $G$,
denote by $G^\sigma$ matrix $(G^\sigma)^i_j = G^{\sigma(i)}_j$.
Let ${\cal K}$ be the quotient skewfield of 
${\cal G} = {\cal NGL}(n,\genfd)$.
We introduce 
the following matrices over ${\cal K}$ and subalgebras of ${\cal K}$:
\begin{itemize}
\item ${}_{r}T$ is the copy of the matrix of generators 
$T \in M_n(C)$ in the image of $M_n(C_r)$ in $M_n(H(C))$ i.e.
$({}_{r}T)^i_j = (t^i_j)_r$. By permuting rows, 
we obtain also ${}_{r}T^\sigma$.
\item $i_\sigma : {\cal G} \rightarrow {\cal W}_{\sigma}$ 
is the Cohn invertive localization of ${\cal G}$
which is universal with respect to the inversion of the
multiplicative set of matrices generated 
by all the lower right square (=principal) submatrices of ${}_{r}T^\sigma$
for all $r$ even and all upper left square (=coprincipal) 
submatrices of ${}_{r}T^\sigma$ for all $r$ odd. Writing ${\cal G}$ is
an inductive limit of Cohn localizations of ${\cal NM}_n$,
and representing localization $i_\sigma$ by inductive limit too,
the same argument as in the proof of Proposition 6, proves that
$i$ is 1-1 and hence ${\cal W}_{\sigma}$ can be viewed as
subalgebra in ${\cal K}$.
\item The Gauss decomposition of matrix 
$T = UA$ in ${\cal W}_{\rm id}$ exists and is unique in
appropriate Cohn localization of ${\cal NM}_n$.
The localization $i_\id$ is 1-1, hence this decomposition
implies decompositions ${}_{r}T = U_r A_r$ if $r$ is even and
${}_{r}T = A_r U_r$ if $r$ is odd, as matrices over 
${\cal W}_{\rm id}$. Here $U_r = ((u_r)^i_j)$ is an upper unitriangular
and $A_r= ((a_r)^i_j)$ a lower triangular matrix. 
If $r$ is ommitted, we mean $r = 0$.
This straightforwardly generalizes to $\sigma \neq\id$. 
\item Let ${\cal U}_\sigma$ be the subalgebra of ${\cal W}_{\sigma}$
generated by the entries of $U_{r,\sigma}$ for all $r$. Let ${\cal A}_r$ be
a subalgebra of ${\cal W}_\sigma$ generated by the entries
of all $A_{r,\sigma}$.
\end{itemize}

\th{Theorem}{7.}{ 1. ${\cal W}_{\sigma}$ is 
generated by ${\cal U}_\sigma$ and ${\cal A}_\sigma$ and 
the natural map ${\cal W}_{\sigma} \rightarrow {\cal K}$ is 1-1.

2. Cohn localization map 
$i_\sigma : {\cal G} \rightarrow {\cal W}_{\sigma}$ 
is 1-1 and ${\cal NB}(n,\genfd)$-coaction compatible.

3. $b^i_j\mapsto (a_{\sigma})^i_j$ extends to
a 1-1 homomorphism of algebras $\gamma_\sigma$.

4. $\gamma_\sigma$ is a ${\cal NB}(n,\genfd)$-comodule map.

5. ${\cal U}_\sigma$ is a strict
subalgebra of the algebra of localized coinvariants
in ${\cal W}_{\sigma}$.
} 

\th{Conjecture.}{}{ The set of $n!$ Cohn localizations 
$i_\sigma : {\cal G} \rightarrow {\cal W}_\sigma$ covers ${\cal G}$.}

\remar{Remarks.}{ A) The quantum analogue (Theorems 9-11) 
does not follow from Theorem 7.
Namely, quantum subgroups of ${\cal NGL}(n,\genfd)$ are obtained
by taking big Hopf ideals, which are not {\it apriori} compatible
with all the structure contained in the theorem. Quantum case
is much harder, particularly for $\sigma \neq \id$. See below.

B) Bruhat decomposition for matrices might have some geometrical
role here, but we did not consider it. 
Algebraically, existence of Bruhat decomposition for matrices
over any skewfield is an easy classical fact~\cite{Cohn:alg3}.
}
Proof~\cite{Skoda:loc-coinv2}
of Theorem 7 (sketch): Perform left row operations
following familiar Gauss elimination procedure, and observe that,
at each step, we can extend the coaction to all the inverses which appear
in the localization corresponding to that step. 
We introduce block matrix
$T[i] = \left(\begin{array}{cc} a &b \\ c& d \end{array}\right)$
where $i\times i$ is the size of submatrix $d$. At the
begining of the induction procedure, $d$ is size $1 \times 1$,
and $T[i] = T$ when we forget the partitioning. 
At each step, using left row operations, change $T[i]$ into 
$A[i] = \left(\begin{array}{cc} a-bd^{-1}c &0 \\ c& d \end{array}\right)$.
Performing these left row operations amounts to the left multiplication
by an upper unitriangular matrix $U[i]$.
Block-entry in the left upper corner 
is simply the $(1,1)$-quasideterminant
of partitioned matrix $T[i]$ (observe, using heredity, which
quasiminors of underlying $n\times n$ matrix
do not change in this procedure). Only matrix $d$
has been inverted at current step. The algebra generated by (genuine,
not block) entries of $A[i]$, contains 
all entries of $d^{-1}$ (use invertibility of $T$),  
and it is equal to the Cohn localization of the algebra 
at submatrix $d$ of the algebra at the previous step.
If $P$ is the basis matrix for noncommutative Borel subgroup ${\cal NB}_n$
(or, in strengthened version of this argument, for the appropriate
parabolic ${\cal NP}_n(i)$), then we write $P$ as a partitoned matrix
(of the same partition type as $T[i]$) 
$P =  \left(\begin{array}{cc} {\frak a} & 0 \\ {\frak c}& {\frak d} 
\end{array}\right)$. Then $\rho(d) = d \otimes {\frak d}$,
hence it is invertible in the localization tensored with ${\cal NB}_n$,
because ${\frak d}$ is invertible in ${\cal NB}_n$, by the simple fact
that a triangular matrix is invertible iff each diagonal entry is such.
Coaction therefore extends to the localization. Moreover,
\[\rho(a-bd^{-1}c) = a \otimes {\frak a} + b \otimes {\frak c}
- (b \otimes {\frak d} )( d^{-1} \otimes {\frak d}^{-1})
(c \otimes {\frak a} + d \otimes {\frak c}) =
(a-bd^{-1}c) \otimes {\frak a}.
\]
Hence $\rho(A[i]) = A[i]\otimes P$ (matrix multiplication with $\otimes$
between the entries).  
This coaction is identical to coaction on matrix $T[i]$,
hence the correspondence 
$\gamma^{i,0} : (T[i])^i_j \rightarrow (A[i])^i_j$ extends to
${\cal NB}(n,\genfd)$-comodule map on the $\genfd$-span of 
set of entries of $T[i]$.
We start with free generators, and gradually impose
relations, hence $\gamma^{i,0}$ extends to 
${\cal NB}(n,\genfd)$-comodule algebra map $\gamma^i$. 
Repartition $A[i]$, making $d$ of size $i+1$, and call
the obtained matrix $T[i+1]$ and continue procedure. Notice
that the composition of the obtained maps $\gamma^i$ 
has the same kernel as the map defining Borel subgroup 
(one inclusion follows from the procedure
and another from generalities on smash products).
These arguments essentially show 1-4. To prove 5. just use the
Gauss decomposition formula for generators of ${\cal U}_\sigma$  as a 
fraction of two flag quasideterminants $u_1 u_2^{-1}$ and observe that
$\rho_\sigma(u_1) = u_1 \otimes b$ and $\rho_\sigma(u_2) = u_2 \otimes b$
for the same $b \in {\cal NB}(n,\genfd)$, and conclude that
$u_1 u_2^{-1}$ is a coinvariant. $\Box$

\section*{14. Quantum matrix groups.}
~\cite{ParshallWang,ManinQGNG,Majid,KlymikSchmud,Hay:qmultilinear,RTF}
Let $q \in \genfd, q \neq 0$. 
The quantum matrix bialgebra
${\cal M}_q(n,\genfd) = {\cal O}(M_q(n,\genfd))$
is the free matrix bialgebra ${\cal NM}(n,\genfd)$ with basis 
$T = (t^\alpha_\beta)$ modulo the smallest 
biideal $I$ such that the following relations hold in quotient:
\[\begin{array}{lcl}
\alpha = \beta, \,\,\gamma < \delta \,\mbox{ (same row) }
&\,\,\,\,\,&
 t^\alpha_\gamma t^\alpha_\delta = q t^\alpha_\delta t^\alpha_\gamma 
\\
\alpha < \beta, \,\,\gamma = \delta \,\mbox{ (same column) }
&&
 t^\alpha_\gamma t^\beta_\gamma = q t^\beta_\gamma t^\alpha_\gamma
\\
\alpha \neq \beta\,\,\mbox{ and }\,\,\gamma \neq \delta
&&
[t^\alpha_\gamma, t^\beta_\delta] = 
 (q - q^{-1}) t^\beta_\gamma t^\alpha_\delta 
(\theta(\delta > \gamma) - \theta (\alpha > \beta)) 
\end{array}\]
where  $\theta({\rm true}) = 1$, $\theta({\rm false}) = 0$,
and $[,]$ stands for the ordinary commutator.

${\cal M}_q(n,\genfd)$ is a domain.
There are other versions, including the {\it multiparametric} case
${\cal M}_{P,Q}(n,\genfd)$, 
cf.~\cite{Manin:multi,strueber,Skoda:thes}.
Some of the results below generalize to ${\cal M}_{P,Q}(n,\genfd)$.

The {\bf quantum determinant} $D \in {\cal M}_q(n,\genfd)$ 
is defined by any of the formulas
\[ D = \sum_{\sigma \in \Sigma(n)} (-q)^{l(\sigma)-l(\tau)} 
t^{\tau(1)}_{\sigma(1)} t^{\tau(2)}_{\sigma(2)} 
\cdots t^{\tau(n)}_{\sigma(n)}
= \sum_{\sigma \in \Sigma(n)} (-q)^{l(\sigma)-l(\tau)} 
t^{\sigma(1)}_{\tau(1)} 
t^{\sigma(2)}_{\tau(2)} \cdots t^{\sigma(n)}_{\tau(n)},\]
where $\Sigma(n)$ is the permutation group on $n$ letters
and $l$ standard length function on $\Sigma(n)$ and
$\tau \in \Sigma(n)$ is a fixed permutation (say identity), 
the choice being irrelevant.
 
$D$ is a central element in ${\cal M}_q(n,\genfd)$. 
$D$ can alternatively be defined using canonical coactions on
quantum exterior algebra~\cite{ParshallWang,Skoda:thes}.

Selecting $m$ rows and columns out of $n$
clearly singles out a subalgebra in ${\cal M}_q(n,\genfd)$ isomorphic
to ${\cal M}_q(m,\genfd)$. {\bf Quantum minors} are the
corresponding quantum determinants $D^K_L$ for so selected
row and column $m$-multilabels $K$ and $L$. Quantum minors
do not necessarily commute with elements which do not belong
to the subalgebra ${\cal M}_q(m,\genfd)$ generated by
selected entries of $T^K_L$. Numerous authors found,
often independently, large number of useful commutation
relations, identities and computational
principles involving quantum minors, 
cf. e.g.~\cite{ParshallWang,Hay:qmultilinear,KrobLec,Skoda:thes}.
The starting point is
\vskip .03in
{\bf Laplace expansion.}
{\it For any pair of ordered 
$m$-multiindices $K,L$ (labels: $1$ to $n$),
\[\begin{array}{lcll}
\delta^K_L D & =& 
\sum_J (-q)^{J-L} D^K_J D^{\hat{L}}_{\hat{J}} 
&= \sum_J (-q)^{J-L} D^J_K  D^{\hat{J}}_{\hat{L}}\\
&=& \sum_J (-q)^{L-J} D^{\hat{L}}_{\hat{J}}D^K_J
&= \sum_J (-q)^{L-J} D^{\hat{J}}_{\hat{L}} D^J_K 
\end{array}\]
where $J$ runs over $m$-multiindices;
the name of a multiindex in exponent denotes the sum 
of labels, and $\,{}\hat{}$ denotes 
the ordered complement with respect to the set 
of labels $\{1,\ldots,n\}$. 
}
\vskip .03in
Define ${\cal SL}_q(n,\genfd) = {\cal O}(SL_q(n,\genfd)) = 
{\cal M}_q(n,\genfd)/L$ where
$L$ is the biideal generated by $(D-1)$. Similarily define
${\cal GL}_q(n,\genfd) = {\cal O}(GL_q(n,\genfd))$ as the 
localization of ${\cal M}_q(n,\genfd)$ at central element $D$. 
In multiparametric ${\cal GL}_{P,Q}(n,\genfd)$ is more interesting
than ${\cal SL}_{P,Q}(n,\genfd)$ because 
the latter typically degenerates.
${\cal SL}_q(n,\genfd)$ and ${\cal GL}_q(n,\genfd)$
are Hopf algebras. The following formula for the antipode
on the generators is forced by Laplace expansion formulas:
\[ St^i_j = (-q)^{i-j} D^{-1} D^{\hat{j}}_{\hat{i}}.\]
Antipode gives the inverse of matrix $T$ in matrix ring over
${\cal GL}_q(n,\genfd)$ by $St^i_j = (T^{-1})^i_j$ and this, in turn,
gives formulas for all the {\it quasideterminants} of $T$.
Namely, $|G|^{-1}_{ij} = (G^{-1})_{ji}$ whenever both
sides are defined. Hence
\[ 
|T|_{ij} = (S(t^j_i))^{-1} = 
(-q)^{j-i} D (D^{\hat{i}}_{\hat{j}})^{-1}
\]
Thus one can use all we know about quasideterminants 
to study quantum determinants, and in particular one can write down
formulas for Gauss decomposition of matrix $T$ and
matrices obtained from $T$ by permutation of rows in
terms of quantum determinants.

\th{Theorem}{8.~\cite{Skoda:thes}}
{Every set of quantum minors
multiplicatively generates a 2-sided Ore set in ${\cal M}_q(n,\genfd)$.
This Ore localization is isomorphic to Cohn localization 
at the set of corresponding submatrices of $T$.
The localization map is 1-1. 
Analogue holds for ${\cal GL}_q(n,\genfd)$, ${\cal SL}_q(n,\genfd)$
and for
(strong) multiparametric deformations ${\cal M}_{P,Q}(n,\genfd)$
and ${\cal GL}_{P,Q}(n,\genfd)$.}

Proof is easy when, for each of the quantum minors involved, 
the selected rows are adjacent (without gaps)
and the selected columns of $T$ are adjacent too. General
case (with gaps) has been proved~\cite{Skoda:thes} 
using several reductions to special cases
and an induction using rather nontrivial commutation relations
involving quantum minors. ``Strong'' denotes a usually imposed
condition on deformation parameters~\cite{artinscheltertate,
demidov:rev, Skoda:thes}.

From now on, ${\cal G}$ will denote either ${\cal GL}_q(n,\genfd)$
or ${\cal SL}_q(n,\genfd)$.

\vskip .005in
{\bf Quantum Borel subgroup} ${\cal B} = {\cal B}_q(n,\genfd)$
is the quotient of ${\cal G}$ by the biideal $I$ generated by $t^i_j$
with $i<j$. $I$ is a Hopf ideal, hence ${\cal B}$ is a Hopf algebra.

\vskip .01in
Let $I$ be a subset $I \subset \{1,\ldots,n-1\}$ and
$I_1, \ldots, I_m$ be its connected components 
ordered in a such a way that
$i< j$ implies every element of $I_i$ 
is smaller than every element of $I_j$.
Denote $I_k^+ = I_k  \cup \{i+1 : i \in I_k \}$ 
for all $1 \leq k \leq m$.
Then the ideal $J_I$ generated by $t^i_j$ where $i < j$ and
$(i,j) \notin \cup_k I^+_k \times I^+_k$ is a Hopf ideal.
If $I$ is nonempty, the quotient Hopf algebra ${\cal G}/J_I$ is called
the {\bf quantum parabolic subgroup} ${\cal P}^I$.
For $I = \emptyset$ we obtain ${\cal P}^I = {\cal B}$ 
and for $I = \{1,\ldots,n-1\}$ we obtain ${\cal P}^I = {\cal G}$.
Notice that $I^+_k \times I^+_k$ is the set of labels of all boxes
in a rectangle with two of the corners at main diagonal.
Hence there is a natural partition of $T$ into 
square submatrices
such that the anihilated ideal is generated by the
elements in the strictly upper block triangular part of $T$.
The dimensions of the blocks are determined by 
the combinatorics of set $I$.
For a given ring $R$, equip the matrices in the matrix ring $M_n(R)$
with that same block partition, and let $M_n^I(R)$
be the resulting ring of block-matrices. 
The identity map surely gives the
ring isomorphism between $M_n(R)$ and $M_n^I(R)$, but the
notions like upper triangular, unidiagonal etc. differ.

\vskip .01in
For $I = \emptyset$ we obtain ${\cal P} = {\cal B}$ 
and for $I = \{1,\ldots,n-1\}$ we obtain ${\cal P} = {\cal G}$.

\vskip .01in
If $I' \subset I$ then $J_I \subset J_{I'}$, hence
there is a natural Hopf algebra surjection 
$\pi_{I,I'} : {\cal P}^I\rightarrow {\cal P}^{I'}$
generalizing $\pi : {\cal G} \rightarrow {\cal B}$
and $\pi_I : {\cal G} \rightarrow {\cal P}^I$.

\vskip .01in
For a permutation $\sigma \in \Sigma(n)$, let $w_\sigma$
be the corresponding permutation matrix,
i.e. for every matrix $G$ we have 
$(w_\sigma^{-1} G)^i_j = g^{\sigma(i)}_j$.
Consider now the noncommutative Gauss decomposition in 
$M_n^{I'}({\cal K}({\cal P}^I))$, where ${\cal K}({\cal P}^I)$
is the quotient skewfield of ${\cal P}^I$:
\begin{equation}\label{eq:par-Gauss}
\pi_I(T) = w_\sigma U_\sigma A_\sigma 
\end{equation}
where $U_\sigma$ is upper unitriangular 
and $A_\sigma$ lower triangular {\it in
$M_n^I({\cal K}({\cal P}^I))$} i.e. in the sense of $I'$-blocks.

Let $(u_\sigma)^i_j$ be the $(i,j)$-th entry of $U_\sigma$,
but now taken as a matrix without block partition
and similarily for $(a_\sigma)^i_j$.

\vskip .004in
For a subset $X$ in some algebra, 
denote by $\langle X \rangle$ the subalgebra generated by $X$. 
Denote 
\begin{equation}\label{eq:calu}\begin{array}{l}
{\cal U}_\sigma = 
\langle (u_\sigma)^i_j \rangle \subset {\cal K}({\cal P}^I),\\
{\cal A}_\sigma = \langle (a_\sigma)^i_j \rangle 
\subset {\cal K}({\cal P}^I),\\
{\cal W}_\sigma = \langle {\cal U}_\sigma , {\cal A}_\sigma \rangle
\subset {\cal K}({\cal P}^I).
\end{array}\end{equation} 

Permutation group $\Sigma(n)$ plays here a role of the Weyl group
of $SL_n$. If we consider only the case of Borel subgroup, 
then the different permutations in $\Sigma(n)$ give different
subalgebras ${\cal U}_\sigma$ in ${\cal G}$. If we consider other parabolics, 
those permutations for which the notion of
being upper $I'$-block triangular agrees, lead only to permutations
among the generators $(u_\sigma)^i_j$ and hence yield 
{\it identical} subalgebras in ${\cal G}$.
Thus for general parabolics
one may consider the relative Weyl group to avoid repetitions.

\vskip .01in
Equation~(\ref{eq:par-Gauss}) can be solved in terms of
quasideterminants, applying formulas for 
noncommutative Gauss decomposition
to $G = w_\sigma^{-1} \pi_I(T)$ considered as a
matrix in $M_n^{I'}({\cal K}({\cal P}^I))$. The general
case is somewhat cumbersome. The simplest, Borel case, as
well as the decomposition where the Hopf ideals are 
generated by all entries in upper triangular part of
first several rows (``fine Gauss decomposition''),
were treated in some detail in author's thesis.
We are now going to formulate corresponding structure theorems
for the Borel case (i.e. analogue of $G/B$)
which give data for a locally trivial quantum principal
fibration in the sense developed in earlier sections
of this paper. 

\th{Theorem}{9.~\cite{Skoda:thes,Skoda:loc-coinv2}}
{Let $S_\sigma$ be the Ore set in ${\cal G} = {\cal GL}_q(n)$
or ${\cal G}={\cal SL}_q(n)$ generated by
all principal (=right lower corner) quantum minors of $w_\sigma^{-1}T$.

(i) The natural coaction 
\[ \rho_{\cal B} : ({\rm id} \otimes \pi)\circ \Delta_{\cal G} : {\cal G} 
\rightarrow {\cal G} \otimes {\cal B} \]
extends to localization making it a ${\cal B}$-comodule algebra
i.e. $S$ is $\rho_{\cal B}$-compatible.

(ii) The quantum Gauss decomposition
\[ T = w_\sigma U_\sigma A_\sigma \]
where $U_\sigma$ is upper unitriangular and $A_\sigma$ lower triangular
has a unique solution for $U_\sigma$ and $A_\sigma$ in matrices over
localized ring $S_\sigma^{-1} {\cal G}$. 

(iii) The natural map of  $S_\sigma^{-1} {\cal G}$ into quotient
skewfield identifies $S_\sigma^{-1} {\cal G}$ with ${\cal W}_\sigma$.
}
\th{Theorem}{10.~\cite{Skoda:thes,Skoda:loc-coinv2}}
{Let $b^i_j = \pi(t^i_j)$ be $(i,j)$-th
generator of quantum Borel. Correspondence 
$b^i_j\rightarrow (a_\sigma)^i_j$ extends to a homomorphism
of ${\cal B}$-comodule algebras
}
\[ \gamma_\sigma : {\cal B} \rightarrow S^{-1}_\sigma {\cal G}. \]
This theorem is relatively easy  
to check for $\sigma = {\rm id}$
i.e. for main cell.
It is harder to show (algebra homomorphism part) for other $\sigma$,
as it amounts to show that $(a_\sigma)^i_j$, which are, up to
constants, ratios of quantum minors, satisfy the relations
which $b^i_j$ do.
There is however a principle (``included row exchange
principle''~\cite{Skoda:thes}) 
which enables to transfer certain identities among
the quantum minors to identities with rows 
sistematically permuted. 
Using this principle, Theorem 10 reduces to the case $\sigma = {\rm id}$.
In the case $\sigma = {\rm id}$ there is an alternative proof. 
Namely, one can use row operations starting with matrix $T$, performing
a version of Gauss elimination procedure which ends with
triangular matrix $((a_{\rm id})^i_j)$. 
One simply shows that at every step
required identities among the entries 
are preserved modulo ideal ${\rm ker}\,\pi$.
However, this ``natural'' argument does not apply to
the case $\sigma \neq {\rm id}$. Namely, in that case, we start
with matrix $w_\sigma^{-1}T$ instead. Its entries
do {\it not} satisfy the required identities modulo ${\rm ker}\,\pi$.
The row operations in that case change the identities too. 
In fact, the situation gets improved, and we
{\it end} with the correct identities! It works but a better
explanation ought to be found.
\th{Theorem}{11.~\cite{Skoda:thes,Skoda:loc-coinv2}} {
(i) ${\cal U}_\sigma \subset (S^{-1}_\sigma{\cal G})^{{\rm co}{\cal B}}$.

(ii) If $\sigma = {\rm id}$ or
if $\sigma = (n\ldots 2\,1)$ then
${\cal U}_\sigma = (S^{-1}_\sigma{\cal G})^{{\rm co}{\cal B}}$.
}
Assertion (i) follows by a direct and simple computation
after expressing $(u_\sigma)^i_j$ in terms of quantum minors.
By Theorem 9 (iii) we see that
${\cal A} =\gamma_\sigma({\cal B})$ 
and ${\cal U}_\sigma$ generate $S^{-1}_\sigma{\cal G}$. Hence (i)
imply by standard arguments
(cf. our discussion of higher smash products above).
that in order to prove (ii) it is sufficient to show
that ${\cal U}_\sigma$ is invariant with respect to
${\cal B}$-action  
\[ 
b \rhd u = \sum \gamma_\sigma(b_{(1)}) u \gamma_\sigma(Sb_{(2)}).
\]
Furthermore, this is enough to show on generators i.e. that
\[ b^i_j \rhd (u_\sigma)^k_l = \sum_{i\geq s\geq j} 
\gamma_\sigma(b^i_s) (u_\sigma)^k_l \gamma_\sigma(Sb^s_j)
 \in {\cal U}_\sigma.\]
We know that $a^i_j =\sigma(b^i_j)$;
we can also express $a^i_j$ in terms of quasideterminants 
(formula for Gauss decomposition!) and then reexpress
in terms of quantum determinants. Here
$Sb^l_j\equiv S_{\cal B}(b^l_j)$ is equal to
$\pi(S_{\cal G}t^l_j)$ and hence, up to a scalar, to the projection
in ${\cal B}$ of a certain quantum (cofactor) minor.
Thus $\gamma_\sigma(Sb^s_j)$ is also a quantum minor
expression but in algebra ${\cal A}_\sigma$. 
If they would commute, say up to a constant factor,
$(a_\sigma)^i_s$ with $(u_\sigma)^k_l$, we would obtain 
a scalar factor times 
\[ (u_\sigma)^k_l
\sum_s \gamma_\sigma(b^i_s) \gamma_\sigma(Sb^s_j) = \delta^i_j u^k_l, \]
where $\delta^i_j$ is Kronecker delta. 
In other words, if $(a_\sigma)^i_s$ and $(u_\sigma)^k_l$ 
commute up to a number depending only on $i,j,k,l,$ 
we obtain an answer proportional to
one obtained from a trivial action $b.u = \epsilon(b)u$.
Natural strategy is hence to try commuting $(a_\sigma)^i_s$ 
(or alternatively, the other, antipode part) 
with $(u_\sigma)^k_l$. The result depends on relative
positions of indices $i,j,k,l$ and it implies Theorem 11 (ii),
what has been checked case by case in author's thesis. 

Assertion (ii) is not a sharp result: we need
the exhaustive list of all 
$n$ and $\sigma \in \Sigma(n)$
when the statement holds. Case by case proof is brute force. 
For general $\sigma$ extra summands appear 
which seem not to live in ${\cal U}_\sigma$.
This phenomenon needs further clarification.

\section*{15. Applications. Comparison with other approaches.}
The concept of a {\it family of Perelomov coherent states}
can be generalized for Hopf algebras~\cite{Skoda:thes,Skoda:coh-states}
in the framework of this paper. Earlier a generalization
We believe that in the special case
of quantum groups the coherent states of~\cite{JurSt:coh}
essentially coincide with our construction.
Under rather general assumptions
a Hopf algebraic analogue of 
a classical resolution of unity by coherent states
has been proved by author in 1999~\cite{Skoda:coh-states}.
\vskip .01in
In the approaches to quantum principal bundles without localization,
noncommutative analogues of the differential calculus and 
of the connections on fibre bundles
were considered in many earlier works, 
e.g.~\cite{BrzMaj:quantumgauge, PflaumSchau, Hajac:strong}. 
{\sc V.~Lunts} and {\sc A.~Rosenberg}~\cite{ LuntsRos1, LuntsRosMP1}
properly extended the {\sc Grothendieck's} definition~\cite{EGAIVend} 
of the rings of regular differential operators 
to noncommutative rings
and considered extensibility of Hopf module algebra actions 
given by the regular differential operators. 
For a left Ore set $S$ in a domain $R$, a derivation $d: R\rightarrow R$ 
always uniquely extends~\cite{Dixmier:env,LuntsRos1} to a derivation 
$d_S : S^{-1}R \rightarrow S^{-1}R$.
Proving that fact is a 
good exercise~\cite{Skoda:loc-coinv1} for a newcomer to Ore localizations.
If a derivation $d : R \rightarrow \Omega^1(R)$
takes value in a $R$-$R$-bimodule $\Omega^1(R)$ such that $(d,\Omega^1(R))$ 
define a 1st order differential calculus, then 
the same proof produces a unique extension 
$d_S : S^{-1}R \rightarrow S^{-1}R \otimes_R \Omega^1(R)$
such that $(d_S, S^{-1}R \otimes_R \Omega^1(R))$ is again
a 1st order differential calculus,
provided an additional condition is satisfied, which we
may call the {\it differential left Ore condition}. Namely, 
$\forall t \in S$, $\forall r \in R$, $\exists s \in S$ and
$\exists \omega \in \Omega^1(R)$, such that  $sdr = \omega t$.

\vskip .007in
We are interested in gluing locally defined 
differential calculi over the quotient spaces, and in the comparison
of the connections defined in quantum vector bundles
over different local charts. For that purpose,
a version of gauge transformations of~\cite{BrzMaj:quantumgauge}
will play role.

\vskip .007in
For a rather different concept of local triviality of noncommutative
principal bundles see~\cite{HajMSz} and references therein.
An early notion of locally trivial vector bundle
using Gabriel localizations appeared in~\cite{Rosen:88}(Appendix 2).
Sheaf-theoretic ideas and localization were present much earlier
in noncommutative geometry, namely in the study of sheaves over
noncommutative spectra~\cite{MurvOyst,vOyst:refl,GolRayOyst}.

\vskip .02in
If $q$ is a primitive root of unity, ${\cal SL}_q(n,{\Bbb C})$ has
a large center. This enables an alternative approach~\cite{semiquantum}
to ${\cal SL}_q(n,{\Bbb C})$
using small noncommutative sheaves over ordinary $SL(n,{\Bbb C})$.
In that case the centers of our localized algebras ${\cal W}_\sigma$
are also big, what enables detailed comparison of the two 
approaches. 

\vskip .03in
Quantum flag varieties can be defined and studied using 
the representation theory of  quantized enveloping 
algebras~\cite{Soib:flag,Soib:book1,Laksh-Resh1,LuntsRos2,Joseph}
or by an Ansatz exploring quantum minors~\cite{TaftTowber}. 
In these approaches, the quantum flag varieties are
described by a noncommutative graded ring ${\cal F}_q$, or by an appropriate 
quotient of a category of (multi)graded ${\cal F}_q$-modules. 
In classical limit $q = 1$, 
this ring can be viewed as the ring of homogeneous functions
on projective variety $G/B$, or as a ring of regular functions
on the quasiaffine variety $G/U$ where $U$ is 
the unipotent radical of $B$.
For type A, ring ${\cal F}_q$ is generated by quantum flag minors.
In our conventions, quantum flag minors are those
quantum $m$-minors which always include all columns from $n-m+1$ to $n$.
The algebra of the main unshifted cell in our construction can
be obtained by taking the zero graded part $(\tilde{S}^{-1}{\cal F}_q)_0$
of a localization of algebra 
${\cal F}_q$ at the homogeneous Ore set $\tilde{S}$ 
generated by the principal flag minors. 
Several authors introduced the quantum flag and Grassmann varieties,
working only 
on the main cell~\cite{Sazdjian:1995yg,JurSt:coh,Stov:Grassmann}.

We note that the Weyl group shifts ${\cal W}_\sigma$ 
are pairwise manifestly isomorphic for $q=1$ but,
in general, nonisomorphic for $q\neq 1$. It would be
interesting, to complete the picture, to find out
whether they are Morita equivalent.

\vskip .015in
Quantum partial flag varieties, and Grassmanians in particular,
can be approached by the methods of this paper. Simple refinements of 
the quantum Gauss decomposition play role there (as parabolics correspond
to block triangular matrices). This enables the induction in stages
and provides the examples of quantum fibre bundles between the different
partial flag varieties corresponding to the inclusions of parabolics
one into another.

\vskip .015in
However, we expect the example ${\cal NGL}_n/{\cal NB}_n$
to be more important,
as the quantum flag varieties seem to be
already pretty tamed by other approaches,
and the fully noncommutative case opens door to a more
unknown area. We also believe in the importance of 
the general foundational questions on quotients, 
coaction--compatible localizations and cohomology.
I hope that a reader of this paper could observe a number of
obvious promising open questions
which the localization approach to quotients opens up.


\section*{Acknowledgements.} 

{\footnotesize
I warmly thank prof. Joel Robbin for his interest and advice
during preparation of my Ph. D. thesis. Several friends,
colleagues and teachers were nice to encourage me 
already early in this project, including 
I.~\v{Z}uti\'{c}, G.~Benkart, B.~Balantekin, 
S.~Witherspoon, A.~Chubukov.  Prof. D.~Passman
scheduled me for a ring theory seminar in Spring 1998,
where first promising
result, Corollary 4 of this paper and $SL_q(2)$
version of Theorems 9-11 were presented.  I thank him
and other people who scheduled or attended my seminars
in Madison, West Lafayette, Zagreb and Bloomington.
Finally, in last about a year, I had a chance to share wider
written, phone and oral communication on the topic
of this work and I thank all who have had patience to
listen or read, and to point out various advice, warnings,
language errors, literature remarks etc. Among them
A.~Ram and C.~Goebel suggested numerous improvements in the
thesis text, and I.~Mori in an after-thesis preprint.
I thank A.~Rosenberg for sharing insights
in a phone conversation, and the same
to A.~Voronov apropos root of unity case.
Y.~Soibelman once suggested a direction
for further applications. Most recently, V.~Lunts
advised me to rework the notion of
quantum coset space via cosimplical techniques (work in progress).
I am thankful to the organizers of the Warsaw conference
``Noncommutative Geometry and Quantum Groups'' (2001) 
for providing a nice workshop atmosphere and reserving a time
slot for my presentation; similarily to the organizers of
``International Hopf Algebra Conference'' at De Paul University
(2002). I am grateful to H.-J.~Schneider for referring me to
the references~\cite{SchMon:quotients,strueber}. 
Thanks to my numerous conversants from the two events,
for their interest and encouragement, particularly to S.~Majid.
I thank to the editor and the referee for the useful suggestions.
Most of the new results in this paper were obtained during
my Ph.D. study at the University of Wisconsin-Madison and
the paper has been completed during my stay at Purdue University.
}


\end{document}